\documentclass[12pt,a4paper,reqno]{amsart}
\usepackage{amssymb}
\usepackage{amscd}
\usepackage{enumerate}
\usepackage{graphicx}
\usepackage{siunitx}
\usepackage{tikz-cd}
\usepackage{color}
\usetikzlibrary{arrows}
\numberwithin{equation}{section}

\usepackage{mathabx}

\usepackage{mathtools}%                  http://www.ctan.org/pkg/mathtools
\usepackage[tableposition=top]{caption}% http://www.ctan.org/pkg/caption
\usepackage{booktabs,dcolumn}%           http://www.ctan.org/pkg/dcolumn + http://www.ctan.org/pkg/booktabs

\DeclareFontFamily{OT1}{rsfs}{}
\DeclareFontShape{OT1}{rsfs}{n}{it}{<-> rsfs10}{}
\DeclareMathAlphabet{\mathscr}{OT1}{rsfs}{n}{it}

\theoremstyle{plain}

\newtheorem{theorem}{Theorem}[section]
\newtheorem{proposition}[theorem]{Proposition}

\theoremstyle{definition}

\newtheorem{remark}[theorem]{Remark}

\newcommand\R{\mathbb{R}}
\newcommand\Z{\mathbb{Z}}

\newcommand\eps{\varepsilon}

\usepackage{algorithm, algorithmic}

\begin{document}

\title[Blowup for supercritical NLW]{Finite time blowup for a supercritical defocusing nonlinear wave system}

\author{Terence Tao}
\address{UCLA Department of Mathematics, Los Angeles, CA 90095-1555.}
\email{tao@math.ucla.edu}

%\email{}


\begin{abstract}  We consider the global regularity problem for defocusing nonlinear wave systems
$$ \Box u = (\nabla_{\R^m} F)(u) $$
on Minkowski spacetime $\R^{1+d}$ with d'Alembertian $\Box := -\partial_t^2 + \sum_{i=1}^d \partial_{x_i}^2$, where the field $u\colon \R^{1+d} \to \R^m$ is vector-valued, and $F\colon \R^m \to \R$ is a smooth potential which is positive and homogeneous of order $p+1$ outside of the unit ball, for some $p >1$.  This generalises the scalar defocusing nonlinear wave (NLW) equation, in which $m=1$ and $F(v) = \frac{1}{p+1} |v|^{p+1}$.  It is well known that in the energy sub-critical and energy-critical cases when $d \leq 2$ or $d \geq 3$ and $p \leq 1+\frac{4}{d-2}$, one has global existence of smooth solutions from arbitrary smooth initial data $u(0), \partial_t u(0)$, at least for dimensions $d \leq 7$.  In this paper we study the supercritical case where $d = 3$ and $p > 5$.  We show that in this case, there exists smooth potential $F$ for some sufficiently large $m$ (in fact we can take $m=40$), positive and homogeneous of order $p+1$ outside of the unit ball, and a smooth choice of initial data $u(0), \partial_t u(0)$ for which the solution develops a finite time singularity.  In fact the solution is discretely self-similar in a backwards light cone.  The basic strategy is to first select the mass and energy densities of $u$, then $u$ itself, and then finally design the potential $F$ in order to solve the required equation.  The Nash embedding theorem is used in the second step, explaining the need to take $m$ relatively large.
\end{abstract}

\maketitle

%%%%%%%%%%%%%%%%%%%%%%%%%%%%%%%%%%%%%%%%%%%%%%%%%

\section{Introduction}

Let $\R^m$ be a Euclidean space, with the usual Euclidean norm $v \mapsto \|v\|_{\R^m}$ and Euclidean inner product $v,w \mapsto \langle v,w \rangle_{\R^m}$.  A function $F\colon \R^m \to \R^n$ is said to be \emph{homogeneous of order $\alpha$} for some real $\alpha$ if we have
\begin{equation}\label{homog}
 F(\lambda v) = \lambda^\alpha F(v)
\end{equation}
for all $\lambda > 0$ and $v \in \R^m$.  In particular, differentiating this at $\lambda=1$ we obtain \emph{Euler's identity}
\begin{equation}\label{euler}
\langle v, (\nabla_{\R^m} F)(v)\rangle_{\R^m} = \alpha F(v)
\end{equation}
where $\nabla_{\R^m}$ denotes the gradient in $\R^m$, assuming of course that the gradient $\nabla_{\R^m} F$ of $F$ exists at $v$.  When $\alpha$ is not an integer, it is not possible for such homogeneous functions to be smooth at the origin unless they are identically zero (this can be seen by performing a Taylor expansion of $F$ around the origin).  To avoid this technical issue, we also introduce the notion of $F$ being \emph{homogeneous of order $\alpha$ outside of the unit ball}, by which we mean that \eqref{homog} holds for $\lambda \geq 1$ and $v \in \R^m$ with $\|v\|_{\R^m} \geq 1$.

Define a \emph{potential} to be a function $F\colon \R^m \to \R$ that is smooth away from the origin; if $F$ is also smooth at the origin, we call it a \emph{smooth potential}.  We say that the potential is \emph{defocusing} if $F$ is positive away from the origin, and \emph{focusing} if $F$ is negative away from the origin.  In this paper we consider nonlinear wave systems of the form
\begin{equation}\label{nlw}
\Box u = (\nabla_{\R^m} F)(u)
\end{equation}
where the unknown field $u\colon \R^{1+d} \to \R^m$ is assumed to be smooth, $\Box = \partial^\alpha \partial_\alpha= -\partial_t^2 + \sum_{i=1}^d \partial_{x_i}^2$ is the d'Alembertian operator on Minkowski spacetime 
$$\R^{1+d} := \{ (t,x_1,\dots,x_d): t,x_1,\dots,x_d \in \R \} = \{ (t,x): t \in \R, x \in \R^d\}$$
with the usual Minkowski metric
$$ \eta_{\alpha \beta} x^\alpha x^\beta = - t^2 + x_1^2 + \dots + x_d^2$$
and the usual Einstein summation, raising, and lowering conventions, $m,d \geq 1$ are integers, and $F\colon \R^m \to \R$ is a smooth potential. This is a Lagrangian field equation, in the sense that \eqref{nlw} is (formally, at least) the Euler-Lagrange equations for the Lagrangian
$$ \int_{\R^{1+d}} \frac{1}{2} \langle \partial^\alpha u, \partial_\alpha u \rangle_{\R^m} + F(u)\ d\eta.$$
We will restrict attention to potentials $F$ which are homogeneous outside of the unit ball of order $p+1$ for some exponent $p>1$.  The well-studied \emph{nonlinear wave equation} (NLW) corresponds to the case when $m=1$ and $F(v) = \frac{|v|^{p+1}}{p+1}$ (for defocusing NLW) or $F(v) = - \frac{|v|^{p+1}}{p+1}$ (for focusing NLW), with the caveat that one needs to restrict $p$ to be an odd integer if one wants these potentials to be smooth.  Later in the paper we will restrict attention to the physical case $d=3$, basically to take advantage of a form of the sharp Huygens' principle.

The natural initial value problem to study here is the Cauchy initial value problem, in which one specifies a smooth initial position $u_0\colon \R^d \to \R^m$ and initial velocity $u_1\colon \R^d \to \R^m$, and asks for a smooth solution $u$ to \eqref{nlw} with $u(0,x) = u_0(x)$ and $\partial_t u(0,x) = u_1(x)$.   Standard energy methods (see e.g. \cite{shatah}) show that for any choice of smooth initial data $u_0, u_1\colon \R^d \to \R^m$, one can construct a solution $u$ to \eqref{nlw} in an open neighbourhood $\Omega$ in $\R^{1+d}$ of the initial time slice $\{ (0,x): x \in \R^d \}$ with this initial data.  Furthermore, either such a solution can be extended to be globally defined in $\R^{1+d}$, or else there is a solution $u$ defined on some open neighbourhood $\Omega$ of $\{ (0,x): x \in \R^d \}$ that ``blow up'' in the sense that they cannot be smoothly continued to some boundary point $(t_*,x_*)$ of $\Omega$.  
The \emph{global regularity problem} for a given choice of potential $F$ asks if the latter situation does not occur, that is to say that for every choice of smooth initial data there is a smooth global solution.  Note that as the equation \eqref{nlw} enjoys finite speed of propagation, there is no need to specify any decay hypotheses on the initial data as this will not affect the answer to the global regularity problem.

For focusing potentials $F$, there are well known blowup examples that show that global regularity fails.  For instance, if $m=1$ and $F\colon \R \to \R$ is given by
\begin{equation}\label{fv-def}
 F(v) := - \frac{2}{(p-1)^2} |v|^{p+1}
\end{equation}
for all $|v| \geq 1$ (and extended arbitrarily in some smooth fashion to the region $|v|<1$ while remaining negative away from the origin), then $F$ is a focusing potential that is homogeneous of order $p+1$ outside of the unit ball, and the function $u: \{ (t,x) \in \R^{1+d}: 0 < t \leq 1\} \to \R$ defined by
\begin{equation}\label{ut}
 u(t,x) := t^{-\frac{2}{p-1}}
\end{equation}
solves \eqref{nlw} but blows up at the boundary $t=0$; applying the time reversal symmetry $(t,x) \mapsto (1-t,x)$, we obtain a counterexample to global regularity for this choice of $F$.  We will thus henceforth restrict attention to defocusing potentials $F$, which excludes ODE-type blowup examples \eqref{ut} in which $u(t,x)$ depends only on $t$. 

The \emph{energy} (or \emph{Hamiltonian})
\begin{equation}\label{energy}
E[u(t)] := \int_{\R^d} \frac{1}{2} \| \partial_t u(t,x)\|_{\R^m}^2 + \frac{1}{2} \| \nabla_x u(t,x)\|_{\R^m \otimes \R^d}^2 + F(u(t,x))\ dx,
\end{equation}
is (formally, at least) conserved by the flow \eqref{nlw}.  A dimensional analysis of this quantity then naturally splits the range of parameters $(d,p)$ into three cases:
\begin{itemize}
\item The \emph{energy-subcritical} case when $d \leq 2$, or when $d \geq 3$ and $p < 1 + \frac{4}{d-2}$.
\item The \emph{energy-critical} case when $d \geq 3$ and $p = 1 + \frac{4}{d-2}$.
\item The \emph{energy-supercritical} case when $d \geq 3$ and $p > 1 + \frac{4}{d-2}$.
\end{itemize}

In the energy-subcritical and energy-critical cases one has global regularity for any defocusing NLW system, at least when $d \leq 7$; see\footnote{Several of these references restrict attention to the scalar NLW or to three spatial dimensions, but the arguments extend without difficulty to the energy-critical NLW systems considered here in the range $3 \leq d \leq 7$.  There are technical difficulties establishing global regularity in extremely high dimension, even when the potential $F$ and all of its derivatives are bounded; see e.g. \cite{wahl}.} \cite{jorgens} for the subcritical case, and \cite{gr1}, \cite{gr2}, \cite{struwe}, \cite{shatah} for the critical case.  These results were also extended to the logarithmically supercritical case (in which the potential $F$ grows faster than the energy-critical potential by a logarithmic factor) in \cite{tao-supercrit}, \cite{roy}. A major ingredient in the proof of global regularity in these cases is the 
conservation of the energy \eqref{energy}, which is non-negative in the defocusing case.  In the energy-critical (and logarithmically supercritical) case, one also takes advantage of Morawetz inequalities such as
\begin{equation}\label{morawetz}
 \int_0^T \int_{\R^d} \frac{F(u(t,x))}{|x|}\ dx dt \leq C E[u(0)] 
\end{equation}
for any time interval $[0,T]$ on which the solution exists.  These bounds can be deduced from the properties of the stress-energy tensor
$$ T_{\alpha \beta} := \langle \partial_\alpha u, \partial_\beta u \rangle - \frac{1}{2} \eta_{\alpha \beta} ( \langle \partial^\gamma u, \partial_\gamma u \rangle_{\R^m} + F(u) )$$
and in particular in the divergence-free nature $\partial^\beta T_{\alpha \beta} = 0$ of this tensor.

It thus remains to address the energy-supercritical case for defocusing smooth potentials $F$.  In this case it is known that the Cauchy problem is ill-posed in various technical senses at low regularities \cite{lebeau}, \cite{cct}, \cite{lebeau-2}, \cite{brenner}, \cite{burq}, \cite{imm}, despite the existence of global weak solutions \cite{segal}, \cite{strauss}, as well as global smooth solutions from sufficiently small initial data \cite{ls} (assuming that $F$ vanishes to sufficiently high order at the origin); see also \cite{zheng} for a partial regularity result.  However, to the authors knowledge, finite time blowup of smooth solutions has not actually been demonstrated for such equations.  The main result of this paper is to establish such a finite time blowup for at least some choices of defocusing potential $F$ and parameters $d,p,m$:

\begin{theorem}[Finite time blowup]\label{main}  Let $d = 3$, let $p > 1 + \frac{4}{d-2}$, and let $m \geq 2 \max(\frac{(d+1)(d+6)}{2}, \frac{(d+1)(d+4)}{2}+5)+2$ be an integer.  Then there exists a defocusing smooth potential $F\colon \R^m \to \R$ that is homogeneous of order $p+1$ outside of the unit ball, and a smooth choice of initial data $u_0, u_1\colon \R^d \to \R^m$, such that there is no global smooth solution $u\colon \R^{1+d} \to \R^m$ to the nonlinear wave system \eqref{nlw}
with initial data $u(0) = u_0, \partial_t u(0) = u_1$.
\end{theorem}

Of course, since $d$ is set equal to $3$, the conditions on $p$ and $m$ reduce to $p > 5$ and $m \geq 40$ respectively.  However, our restriction to the $d=3$ case is largely for technical reasons (basically in order to exploit the strong Huygens principle), and we believe the results should extend to higher values of $d$, with the indicated constraints on $d$ and $p$, though we will not pursue this matter here.  The rather large value of $m$ is due to our use of the Nash embedding theorem (!) at one stage of the argument.  It would of course be greatly desirable to lower the number $m$ of degrees of freedom down to $1$, in order to establish blowup for the scalar defocusing supercritical NLW, but our methods crucially need a large value of $m$ in order to ensure that a certain map from a $1+d$-dimensional space into the sphere $S^{m-1}$ is embedded, which is where the Nash embedding theorem comes in.  Nevertheless, even though Theorem \ref{main} does not directly show that the scalar defocusing supercritical NLW exhibits finite time blowup, it does demonstrate a significant \emph{barrier} to any attempt to prove global regularity for this equation, as such an attempt must necessarily use some special property of the scalar equation that is not shared by the more general system \eqref{nlw}.

We briefly discuss the methods used to prove Theorem \ref{main}.  The singularity constructed is a discretely self-similar blowup in a backwards light cone; see the reduction to Theorem \ref{main-2} below.  In particular, the blowup is ``locally of type II'' in the sense that scale-invariant norms inside the light cone stay bounded, but not ``globally of type II'', as a significant amount of energy (as measured using scale-invariant norms) radiates out of the backwards light cone at all scales.  This is compatible with the results in
\cite{kenig}, \cite{killip}, \cite{killip-2}, which rule out ``global'' type II blowup, but not ``local'' type II blowup.
It would be natural to seek a \emph{continuously} self-similar smooth blowup solution, but it turns out\footnote{On the other hand, it is possible to use perturbative methods to create \emph{rough} solutions to \eqref{nlw} that are continuously self-similar: see \cite{planchon}, \cite{ribaud}.  However, these methods do not seem to be adaptable to generate smooth solutions, and indeed Proposition \ref{nose} suggests that there are strong obstacles in trying to create such an adaptation.  The negative result here also stands in contrast to the situation of high-dimensional wave maps into negatively curved targets, where ODE methods were used in \cite{cst} to construct continuously self-similar blowup examples in seven and higher spatial dimensions.} that these are ruled out; see Proposition \ref{nose} below.  Hence we will not restrict attention to continuously self-similar solutions.  It also turns out to be convenient not to initially restrict attention to spherically symmetric solutions, although we will eventually do so later in the argument.  

Traditionally, one thinks of the potential $F$ as being prescribed in advance, and the field $u$ as the unknown to be solved for.  However, as we have the freedom to \emph{select} $F$ in Theorem \ref{main}, it turns out to be more convenient to prescribe $u$ first, and only then design an $F$ for which the equation \eqref{nlw} is obeyed.  This turns out to be possible as long as the map $\theta \colon (t,x) \mapsto \frac{u(t,x)}{\|u(t,x)\|_{\R^m}}$ has certain non-degeneracy properties, and if the stress-energy tensor $T_{\alpha \beta}$ (which can be defined purely in terms of $u$) is divergence-free; see the reduction to Theorem \ref{main-3} below.  The stress-energy tensor $T_{\alpha \beta}$ (or more precisely, some related fields which we call the \emph{mass density} $M$ and the \emph{energy tensor} $E_{\alpha \beta}$) can be viewed as prescribing the metric geometry of the map $\theta$, and the Nash embedding theorem can then be used to locate a choice of $\theta$ with the desired non-degeneracy properties and the prescribed metric, so long as the fields $M$ and $E_{\alpha \beta}$ obey a number of conditions (one of which relates to the divergence-free nature of the stress-energy tensor, and another to the positive definiteness of the Gram matrix of $u$).  This reduces the problem to a certain ``semidefinite program'' (see Theorem \ref{main-4}), in which one now only needs to specify the fields $M$ and $E_{\alpha \beta}$, rather than the original field $u$ or the potential $F$.

It is at this point (after some additional technical reductions in which certain fields are allowed to degenerate to zero) that it finally becomes convenient to make symmetry reductions, working with fields $M$, $E_{\alpha \beta}$ that are both continuously self-similar and spherically symmetric, and assuming that there are no angular components to the energy tensor.  In three spatial dimensions, this reduces the divergence-free nature of the stress-energy tensor to a single transport equation for the null energy $e_+$ (which, in terms of the original field $u$, is given in polar coordinates by $e_+ = \frac{1}{2} \| (\partial_t + \partial_r) (ru) \|_H^2$), in terms of a certain ``potential energy density'' $V$ (which, in terms of the original data $u$ and $F$, is given by $V = r F(u)$); see Theorem \ref{main-6} for a precise statement.  The strategy is then to solve for these fields $e_+, V$ first, and then choose all the remaining unknown fields in such a way that the remaining requirements of the semidefinite program are satisfied.  This turns out to be possible if the fields $e_+, V$ are chosen to concentrate close to the boundary of the light cone.

The author is supported by NSF grant DMS-1266164 and by a Simons Investigator Award.  The author is also indebted to Sergiu Klainerman for suggesting this question, and Yi Jin Gao for corrections.

\section{Reduction to discretely self-similar solution}\label{reduct}

We begin the proof of Theorem \ref{main}. 

We first observe that from finite speed of propagation and the symmetries of the equation, Theorem \ref{main} follows from the following claim, in which the solution is restricted to a truncated light cone and is discretely self-similar and the potential is now homogeneous everywhere (not just outside of the unit ball), but no longer required to be smooth.  This reduction does not use any of the hypotheses on $m,d,p$.

\begin{theorem}[First reduction]\label{main-2}  Let $d = 3$, let $p > 1 + \frac{4}{d-2}$, and let $m \geq 2 \max(\frac{(d+1)(d+6)}{2}, \frac{(d+1)(d+4)}{2}+5)+2$ be an integer.  Then there exists a defocusing potential $F\colon \R^d \to \R^m$ which is homogeneous of order $p+1$ and a smooth function $u \colon \Gamma_d \to \R^m \backslash \{0\}$ on the light cone $\Gamma_d := \{ (t,x) \in \R^{1+d}: t>0; |x| \leq t \}$ that solves \eqref{nlw} on its domain and is nowhere vanishing, and also discretely self-similar in the sense that there exists $S>0$ such that
\begin{equation}\label{uts}
 u(e^S t, e^S x ) = e^{-\frac{2}{p-1}S} u(t,x)
\end{equation}
for all $(t,x) \in \Gamma_d$.
\end{theorem}

A key point here is that $u$ is smooth all the way up to the boundary of the light cone $\Gamma_d$, rather than merely being smooth in the interior.  The exponent $-\frac{2}{p-1}$ is mandated by dimensional analysis considerations.  It would be natural to consider solutions that are continuously self-similar in the sense that \eqref{uts} holds for \emph{all} $S \in \R$, but as we shall shortly see, it will not be possible to generate such solutions in the three-dimensional defocusing setting.

Let us assume Theorem \ref{main-2} for the moment, and show how it implies Theorem \ref{main}.  Let $F, S, u$ be as in Theorem \ref{main-2}.   Since $u$ is smooth and non-zero on the compact region $\{ (t,x) \in \Gamma_d: e^{-S} \leq t \leq 1 \}$, it is bounded from below in this region.  By replacing $u$ with $Cu$ and $F$ with $v \mapsto C^2 F(v/C)$ for some large constant $C$, we may thus assume that
$$ \| u(t,x) \|_{\R^m} \geq 1$$
whenever $(t,x) \in \Gamma_d$ with $e^{-S} \leq t \leq 1$.  Using the discrete self-similarity property \eqref{uts}, we then have this bound for all $0 < t \leq 1$; in fact we have a lower bound on $\|u(t,x)\|_{\R^m}$ that goes to infinity as $t \to 0$, ensuring in particular that $u$ has no smooth extension to $(0,0)$.

Using a smooth cutoff function, one can find a smooth defocusing potential $\tilde F\colon \R^m \to \R$ that agrees with $F$ in the region $\{ v \in \R^m: \|v\|_{\R^m} \geq 1 \}$.  Then $u$ solves \eqref{nlw} in the truncated light cone $\{ (t,x) \in \R^{1+d}: 0 < t \leq 1; |x| \leq t \}$  with $F$ replaced by $\tilde F$.  Choose smooth initial data $v_0,v_1\colon \R^d \to \R^m$ such that
$$ v_0(x) = u(1,x)$$
and
$$ v_1(x) = - \partial_t u(1,x) $$
for all $|x| \leq 1$ (where we use $|x| := \|x\|_{\R^d}$ to denote the magnitude of $x \in \R^d$); such data exists from standard smooth extension theorems (see e.g. \cite{seeley}) since the functions $u(1,x), \partial_t u(1,x)$ are smooth on the closed ball $\{ x: |x| \leq 1 \}$.  Suppose for contradiction that Theorem \ref{main} failed (with $F$ replaced by $\tilde F$), then we have a global smooth solution $v\colon \R^{1+d} \to \R^m$ to \eqref{nlw} (for $\tilde F$) with initial data $v(0) = v_0, \partial_t v(0) = v_1$.  The function $\tilde u \colon (t,x) \mapsto v(1-t,x)$ is then another global smooth solution to \eqref{nlw} (for $\tilde F$) such that $\tilde u(1,x) = u(1,x)$ and $\partial_t \tilde u(1,x) = \partial_t u(1,x)$ for all $|x| \leq 1$.  Finite speed of propagation (see e.g. \cite[Proposition 3.3]{tao-book}) then shows that $\tilde u$ and $u$ agree in the region $\{ (t,x) \in \R^{1+d}: 0 < t \leq 1; |x| \leq t \}$; as $\tilde u$ is smoothly extendible to $(0,0)$, $u$ is also, giving the desired contradiction.   This concludes the derivation of Theorem \ref{main} from Theorem \ref{main-2}.

It remains to prove Theorem \ref{main-2}.  This will be the focus of the remaining sections of the paper.  For now, let us show why continuously self-similar solutions are not available in the defocusing case, at least for some choices of parameters $d,p$.  The point will be that continuous self-similarity gives a new monotonicity formula for a certain quantity $f(t,r)$ (measuring a sort of ``equipartition of energy'') that can be used to derive a contradiction.

\begin{proposition}[No self-similar defocusing solutions]\label{nose}  Let $d \geq 3$ and $p > 1$ be such that $\frac{d-3}{2} - \frac{2}{p-1} < 0$, let $m$ be a natural number, and let $F\colon \R^m \to \R$ be a defocusing potential that is homogeneous of order $p+1$.  Then there does not exist a smooth solution $u\colon \Gamma_d \to \R^m \backslash \{0\}$ to \eqref{nlw} that is homogeneous of order $-\frac{2}{p-1}$.
\end{proposition}

Note in particular that in the physical case $d=3$, the condition $\frac{d-3}{2} - \frac{2}{p-1} < 0$ is automatic, and so no self-similar defocusing solutions exist in this case.  We do not know if this condition is necessary in the above proposition.

\begin{proof}  Suppose for contradiction that such a $u$ exists.  The equation \eqref{nlw} in polar coordinates $(t,r,\omega)$ reads
$$ -\partial_{tt} u + \partial_{rr} u + \frac{d-1}{r} \partial_r u + \frac{1}{r^2} \Delta_\omega u = (\nabla F)(u) $$
where $\Delta_\omega$ is the Laplace-Beltrami operator on the sphere $S^{d-1}$.  Making the substitution
\begin{equation}\label{phitr}
\phi(t,r,\omega) := r^{(d-1)/2} u(t,r,\omega),
\end{equation}
this becomes
\begin{equation}\label{star}
 -\partial_{tt} \phi + \partial_{rr} \phi - \frac{1}{r^2} \left(-\Delta_\omega + \frac{(d-1)(d-3)}{4}\right) \phi = r^{(d-1)/2} (\nabla F)(r^{-(d-1)/2} \phi) 
\end{equation}
for $r>0$.  

We introduce the scaling vector field $S := t \partial_t + r \partial_r$ and the Lorentz boost $L := r \partial_t + t \partial_r$.  Observe that $L$ and $S$ commute with
\begin{equation}\label{sas}
 - S^2 + L^2 = (t^2-r^2) (-\partial_{tt} + \partial_{rr})
\end{equation}
and thus
$$ - \langle S^2 \phi, L \phi \rangle_{\R^m} + \langle L^2 \phi, L \phi \rangle_{\R^m} = (t^2-r^2) \langle -\partial_{tt} \phi + \partial_{rr} \phi, L\phi \rangle_{\R^m}.$$
As $u$ is assumed homogeneous of order $-\frac{2}{p-1}$, $\phi$ is homogeneous of order $\frac{d-1}{2}-\frac{2}{p-1}$.  
From Euler's identity \eqref{euler} we thus have $\phi$ an eigenfunction of $S$,
$$ S \phi = \left(\frac{d-1}{2}-\frac{2}{p-1}\right) \phi,$$ 
and thus (by the commutativity of $L$ and $S$)
$$ \langle L \phi, S^2 \phi \rangle_{\R^m} = \langle LS \phi, S \phi \rangle_{\R^m} = \frac{1}{2} L \| S \phi \|_{\R^m}^2.$$
We also have
$$ \langle L \phi, L^2 \phi \rangle_{\R^m} = \frac{1}{2} L \| L \phi \|_{\R^m}^2.$$
Putting all of these facts together, we conclude that
$$ L \left( - \frac{1}{2} \| S \phi \|_{\R^m}^2 + \frac{1}{2} \| L \phi \|_{\R^m}^2\right) = (t^2-r^2)  \langle -\partial_{tt} \phi + \partial_{rr} \phi, L\phi \rangle_{\R^m}.$$
A computation similar to \eqref{sas} shows that
$$- \| S \phi \|_{\R^m}^2 + \| L \phi \|_{\R^m} = (t^2-r^2) (- \| \partial_t \phi \|_{\R^m}^2 + \| \partial_r \phi \|_{\R^m}^2 ).$$
Since $t^2-r^2$ is annihilated by $L$, we conclude that
$$ L \left( - \frac{1}{2} \|\partial_t \phi \|_{\R^m}^2 + \frac{1}{2} \| \partial_r \phi \|_{\R^m}^2 \right) = 
\langle -\partial_{tt} \phi + \partial_{rr} \phi, \phi \rangle_{\R^m}.$$
By \eqref{star}, the right-hand side is equal to
$$
\frac{1}{r^2} \langle -\Delta_\omega \phi, L\phi \rangle_{\R^m} + \frac{(d-1)(d-3)}{4 r^2} \langle \phi, L\phi \rangle_{\R^m} + r^{(d-1)/2} \langle (\nabla F)(r^{-(d-1)/2} \phi), L \phi \rangle_{\R^m}.$$
To deal with the angular Laplacian, we integrate over $S^{d-1}$ and then integrate by parts to conclude that
\begin{align*}
 &L \int_{S^{d-1}} ( - \frac{1}{2} \|\partial_t \phi \|_{\R^m}^2 + \frac{1}{2} \| \partial_r \phi \|_{\R^m}^2 ) \ d\omega \\
&\quad = 
\int_{S^{d-1}} 
\frac{1}{2r^2} L \| \nabla_\omega \phi\|_{\R^m \otimes \R^d}^2 + \frac{(d-1)(d-3)}{8 r^2} L \| \phi \|_{\R^m}^2 + r^{(d-1)/2} \langle (\nabla F)(r^{-(d-1)/2} \phi), L \phi \rangle_{\R^m}\ d\omega
\end{align*}
where we use the fact that the Lorentz boost $L$ commutes with angular derivatives, and where $d\omega$ denotes surface measure on $S^{d-1}$.

From the chain and product rules, noting that $Lr = t$, we have
$$ L \phi = r^{(d-1)/2} L( r^{-(d-1)/2} \phi ) + \frac{d-1}{2} \frac{t}{r} \phi$$
and thus (using \eqref{euler})
\begin{align*}
 \langle (\nabla F)(r^{-(d-1)/2} \phi), L \phi \rangle_{\R^m} &= r^{(d-1)/2} \left( L F( r^{-(d-1)/2} \phi ) + \frac{d-1}{2} \frac{t}{r} \langle r^{-(d-1)/2} \phi, (\nabla F)(r^{-(d-1)/2} \phi) \rangle_{\R^m} \right) \\
&= r^{(d-1)/2} \left( L F( r^{-(d-1)/2} \phi ) + \frac{(d-1)(p+1)}{2} \frac{t}{r} F( r^{-(d-1)/2} \phi ) \right).
\end{align*}
Putting all this together, we see that if we introduce the quantity
\begin{align*}
 f(t,r) &:= \int_{S^{d-1}} - \frac{1}{2} \|\partial_t \phi \|_{\R^m}^2 + \frac{1}{2} \| \partial_r \phi \|_{\R^m}^2 \\
&\quad - \frac{1}{2r^2} \| \nabla_\omega \phi\|_{\R^m \otimes \R^d}^2  -  \frac{(d-1)(d-3)}{8 r^2} \| \phi \|_{\R^m}^2 \\
&\quad - r^{d-1} F( r^{-(d-1)/2} \phi )\ d\omega
\end{align*}
then we have the formula
\begin{align*}
Lf &= \int_{S^{d-1}} \frac{t}{r^3} \| \nabla_\omega \phi\|_{\R^m \otimes \R^d}^2 + \frac{(d-1)(d-3)t}{4 r^3} \| \phi \|_{\R^m}^2 \\
&\quad + \frac{(d-1)(p-1)}{2} \frac{t}{r} r^{d-1} F( r^{-(d-1)/2} \phi ) )\ d\omega.
\end{align*}
for any $r>0$. In particular, $f( \cosh y, \sinh y)$ is a strictly function of $y$ for $y>0$, since
$$ \frac{d}{dy} f( \cosh y, \sinh y) = (Lf)(\cosh y, \sinh y) > 0$$
with the strict positivity coming from the defocusing nature of $F$.
On the other hand, when $y \to 0^+$, we see from \eqref{phitr} that all the negative integrands in the definition of $f(\cosh y, \sinh y)$ go to zero, and thus
$$ \lim_{y \to 0^+} f(\cosh y, \sinh y) \ge 0.$$
Combining these two facts, we conclude in particular that 
\begin{equation}\label{true}
\lim_{y \to +\infty} f( \cosh y, \sinh y ) > 0.
\end{equation}
On the other hand, as $\phi$ is homogeneous of order $\frac{d-1}{2}-\frac{2}{p-1}$ and $F$ is homogeneous of order $p+1$, we see that the integrand in the definition of $f(t,r)$ is homogeneous of order $2 (\frac{d-3}{2} - \frac{2}{p-1})$, which is negative by hypothesis.  This implies that $f(\cosh y, \sinh y)$ goes to zero as $y \to +\infty$, contradicting \eqref{true}.
\end{proof}

\section{Eliminating the potential}\label{pot}

We now exploit the freedom to select the defocusing potential $F$ by eliminating it from the equations of motion.  To motivate this elimination, let us temporarily make the \emph{a priori} assumption that we have a solution $u$ to \eqref{nlw} in the light cone $\Gamma_d$ from Theorem \ref{main-2} that is nowhere vanishing.  Taking the inner product of \eqref{nlw} with $u$ and using \eqref{euler} then gives an equation for $F(u)$:
\begin{equation}\label{fum}
F(u) = \frac{1}{p+1} \langle u, \Box u \rangle_{\R^m}.
\end{equation}
In particular, since $F$ is defocusing and $u$ is nowhere vanishing, we have the defocusing property
\begin{equation}\label{pos}
\langle u, \Box u \rangle_{\R^m} > 0
\end{equation}
throughout $\Gamma_d$.  Next, if $\partial_\alpha$ denotes one of the $d+1$ derivative operators $\partial_t, \partial_{x_1}, \dots, \partial_{x_d}$, we have from the chain rule that
$$ \partial_\alpha F(u) = \langle \partial_\alpha u, (\nabla F)(u) \rangle_{\R^m}$$
and hence from \eqref{nlw}, \eqref{fum} we have the equation
\begin{equation}\label{alpha-eq}
\partial_\alpha \langle u, \Box u \rangle_{\R^m} = (p+1) \langle \partial_\alpha u, \Box u \rangle_{\R^m}.
\end{equation}

\begin{remark}  One can rewrite the equation \eqref{alpha-eq} in the more familiar form
$$ \partial^\beta T_{\alpha \beta} = 0$$
where $T_{\alpha\beta}$ is the stress-energy tensor
$$ T_{\alpha \beta} = \langle \partial_\alpha u, \partial_\beta u \rangle_{\R^m} - \eta_{\alpha \beta} \left( \frac{1}{2} \langle \partial^\gamma u, \partial_\gamma u \rangle_{\R^m} + \frac{1}{p+1} \langle u, \Box u \rangle \right).$$
\end{remark}

Now assume that $u$ obeys the discrete self-similarity hypothesis \eqref{uts}.
Let $\theta :=  u / \|  u \|$ denote the direction vector of $ u$, then $\theta$ is smooth map from $\Gamma_d$ to the unit sphere $S^{m-1} := \{ v \in \R^m: \|v\|_{\R^m} = 1 \}$ of $\R^m$.  From the discrete self-similarity \eqref{uts} we see that $\theta$ is invariant under the dilation action of the multiplicative group $e^{S\Z} := \{ e^{nS}: n \in \Z \}$ on $\Gamma_d$.  Thus $\theta$ descends to a smooth map $\tilde \theta\colon \Gamma_d/e^{S\Z} \to S^{m-1}$ on the compact quotient $\Gamma_d/e^{S\Z}$ which is a smooth surface with boundary (diffeomorphic to the product of a $d$-dimensional closed ball and a circle).  Under some non-degeneracy hypotheses on this map, we can now eliminate the potential $F$, reducing Theorem \ref{main-2} to the following claim:

\begin{theorem}[Second reduction]\label{main-3}   Let $d = 3$, let $p > 1 + \frac{4}{d-2}$, and let $m \geq 2 \max(\frac{(d+1)(d+6)}{2}, \frac{(d+1)(d+4)}{2}+5)+2$ be an integer.  Then there exists $S > 0$ and a smooth nowhere vanishing function $u\colon \Gamma_d \to \R^m \backslash \{0\}$ which is discretely self-similar in the sense of \eqref{uts} and obeys the defocusing property \eqref{pos} and the equations \eqref{alpha-eq} throughout $\Gamma_d$.  Furthermore, the map $\tilde \theta\colon \Gamma_d/\lambda^\Z \to S^{m-1}$ defined as above is injective, and immersed in the sense that the $d+1$ derivatives $\partial_\alpha \theta(t,x)$ for $\alpha=0,\dots,d$ are linearly independent in $\R^m$ for each $(t,x) \in \Gamma_d$.
\end{theorem}

Let us assume Theorem \ref{main-3} for now and see how it implies Theorem \ref{main-2}.  As in the previous section, our arguments here will not depend on our hypotheses on $m,d,p$.

Since the map $\tilde \theta\colon \Gamma_d/e^{S\Z} \to S^{m-1}$ is assumed to be injective and immersed, it is a smooth embedding of the set $\Gamma_d/e^{S\Z}$ to $S^{m-1}$, so that $\tilde \theta(\Gamma_d/e^{S\Z}) = \theta(\Gamma_d)$ is a smooth manifold with boundary contained in $S^{m-1}$.  We define a function $F_0 \colon \theta( \Gamma_d ) \to \R$ by the formula
\begin{equation}\label{F0-def}
F_0\left( \frac{u(t,x)}{\|u(t,x)\|_{\R^m}} \right) := \frac{1}{(p+1) \|u(t,x)\|_{\R^m}^{p+1}} \langle u(t,x), \Box u(t,x) \rangle_{\R^m}
\end{equation}
for any $(t,x) \in \Gamma_d$.  As $\theta$ is injective and $u$ is spherically symmetric, nowhere vanishing, and discretely self-similar, one verifies that $F_0$ is well-defined.  As the map $\theta$ is immersed, we also see that $F_0$ is smooth.  From \eqref{pos} we see that $F_0$ is positive on $\theta(\Gamma_d)$.  Intuitively, $F_0$ is going to be our choice for $F$ on the set $\theta(\Gamma_d)$ (this choice is forced upon us by \eqref{fum} and homogeneity).

We define an auxiliary function $T \colon \theta( \Gamma_d) \to \R^m$ by the formula
\begin{equation}\label{T-def}
 T\left( \frac{u(t,x)}{\|u(t,x)\|_{\R^m}} \right) := \frac{1}{\|u(t,x)\|_{\R^m}^{p}} \Box u(t,x) - \frac{1}{\|u(t,x)\|_{\R^m}^{p+2}} \langle u(t,x), \Box u(t,x) \rangle_{\R^m} u(t,x) 
\end{equation}
for all $(t,x) \in \Gamma_d$; geometrically, this is the orthogonal projection of $\frac{1}{\|u\|_{\R^m}^p} \Box u$ to the tangent plane of $S^m$ at $\frac{u}{\|u\|_{\R^m}}$, and will be our choice for the $S^{m-1}$ gradient 
$$(\nabla_{S^{m-1}} F)\left(\frac{u}{\|u\|}\right) = (\nabla_{\R^m} F)\left(\frac{u}{\|u\|_{\R^m}}\right) - \left\langle \frac{u}{\|u\|_{\R^m}}, (\nabla_{\R^m} F)\left(\frac{u}{\|u\|}\right)  \right\rangle_{\R^m} \frac{u}{\|u\|_{\R^m}}$$
of $F$ at $\frac{u}{\|u\|_{\R^m}}$.

As $\theta$ is injective and $u$ is nowhere vanishing and discretely self-similar, one verifies as before that $T$ is well-defined, and from the immersed nature of $\theta$ we see that $T$ is smooth.  Clearly $T(\omega)$ is also orthogonal to $\omega$ for any $\omega \in \theta( \Gamma_d)$.  We also claim that $T$ is an extension of the gradient $\nabla_{\theta( \Gamma_d)} F_0$ of $F_0$ on $\theta( \Gamma_d)$, in the sense that
\begin{equation}\label{gorad}
 \langle v, \nabla_{\theta( \Gamma_d)} F_0(\omega) \rangle_{\R^m} = \langle v, T(\omega) \rangle_{\R^m}
\end{equation}
for any $\omega \in \theta( \Gamma_d)$ and tangent vectors $v \in T_\omega \theta( \Gamma_d)$ to $\theta( \Gamma_d)$ at $\omega$.  To verify \eqref{gorad}, we write $\omega = \frac{u(t,x)}{\|u(t,x)\|_{\R^m}} = \frac{u}{\|u\|}$ for some $(t,x) \in \Gamma_d$; henceforth we suppress the explicit dependence on $(t,x)$ for brevity.  The tangent space to $\theta(\Gamma_d)$ at $\omega$ is spanned by $\partial_\alpha \frac{u}{\|u\|}$ for $\partial_\alpha = \partial_t, \partial_{x_1},\dots,\partial_{x_d}$, so it suffices to show that
$$
\left \langle \partial_\alpha \frac{u}{\|u\|}, \nabla_{\theta( \Gamma_d)} F_0(\omega) \right\rangle_{\R^m} =\left \langle \partial_\alpha \frac{u}{\|u\|}, T(\omega) \right\rangle_{\R^m}
$$
for each $\partial_\alpha$.  But from the chain and product rules and \eqref{F0-def}, \eqref{alpha-eq}, \eqref{T-def} we have
\begin{align*}
 \left\langle \partial_\alpha \frac{u}{\|u\|}, \nabla_{\theta( \Gamma_d)} F_0(\omega) \right\rangle_{\R^m} 
&= \partial_\alpha F_0\left( \frac{u}{\|u\|} \right) \\
&= \frac{1}{p+1} \partial_\alpha \left( \frac{1}{\|u\|_{\R^m}^{p+1}} \langle u, \Box u \rangle_{\R^m} \right) \\
&= -\frac{\langle u, \partial_\alpha u \rangle_{\R^m}}{\|u\|_{\R^m}^{p+3}} \langle u, \Box u \rangle_{\R^m}  + \frac{\partial_\alpha \langle u, \Box u \rangle_{\R^m} }{(p+1)\|u\|_{\R^m}^{p+1}} \\
&= -\frac{\langle u, \partial_\alpha u \rangle_{\R^m}}{\|u\|_{\R^m}^{p+3}} \langle u, \Box u \rangle_{\R^m} + \frac{\langle \partial_\alpha u, \Box u \rangle_{\R^m} }{\|u\|_{\R^m}^{p+1}} \\
&= \left\langle \frac{1}{\|u\|_{\R^m}} \partial_\alpha u, T( \frac{u}{\|u\|_{\R^m}} ) \right\rangle_{\R^m} \\
&= \left\langle \partial_\alpha \frac{u}{\|u\|_{\R^m}}, T\left( \frac{u}{\|u\|_{\R^m}} \right) \right\rangle_{\R^m} 
\end{align*}
as desired, where in the final line comes from the orthogonality of $T\left( \frac{u}{\|u\|_{\R^m}} \right)$ with scalar multiples of $u$.

We now claim that we may find an open neighbourhood $U$ of $\theta(\Gamma_d)$ in $S^{m-1}$ and a smooth extension $F_1 \colon U \to \R$ of $F_0$, with the property that 
\begin{equation}\label{grad} \nabla_{S^{m-1}} F_1(\omega) = T(\omega) 
\end{equation}
for all $\omega \in \theta(\Gamma_d)$.  Indeed, we can define
$$ F_1(\omega + v) := F_0(\omega) + \langle v, T(\omega) \rangle_{\R^m}$$
for all $\omega \in \theta(\Gamma_d)$ and sufficiently small $v \in \R^m$ orthogonal to the tangent space $T_\omega \theta(\Gamma_d/e^{S\Z})$ with $\omega+v \in S^{m-1}$; one can verify that this is well-defined as a smooth extension of $F_0$ to a sufficiently small normal neighbourhood of $\theta(\Gamma_d)$ with the desired gradient property \eqref{grad} (here we use \eqref{gorad} to deal with tangential components of the gradient), and one may smoothly extend this to an open neighbourhood of $\theta(\Gamma_d)$ by Seeley's theorem \cite{seeley}.

Next, if we extend $F_1$ by zero to all of $S^{m-1}$ and define $F_2 \colon S^{m-1} \to \R$ to be the function $F_2 := \psi F_1 + (1-\psi)$ for some smooth function $\psi \colon S^{m-1} \to [0,1]$ supported in $U$ that equals $1$ on a neighbourhood of $\theta(\Gamma_d)$, then $F_2$ is a smooth extension of $F_0$ to $S^{m-1}$ that is strictly positive, and which also obeys \eqref{grad}.  If we then set $F\colon \R^m \to \R$ to be the function
$$ F( \lambda \omega ) := \lambda^{p+1} F_2(\omega) $$
for all $\lambda \geq 0$ and $\omega \in S^{m-1}$, then $F$ is a defocusing potential, homogeneous of order $p+1$, which extends $F_0$, and such that
$$ \nabla_{S^{m-1}} F(\omega) = T(\omega)$$
for $\omega \in \theta(\Gamma_d)$.  By homogeneity \eqref{homog}, the radial derivative $\langle \omega, \nabla_{\R^m} F(\omega)\rangle_{\R^m}$ is $(p+1) F(\omega) = (p+1) F_0(\omega)$ for such $\omega$, and hence for $\omega = \frac{u}{\|u\|}$ by \eqref{T-def}, \eqref{F0-def}
\begin{align*}
\nabla_{\R^m} F(\omega) &= T(\omega) + (p+1) F_0(\omega) \omega \\
&= \frac{1}{\|u\|^p} \Box u - \frac{\langle u,\Box u \rangle}{\|u\|^{p+2}}  u + \frac{p+1}{(p+1)\|u\|^{p+1}} \langle u, \Box u \rangle \frac{u}{\|u\|} \\
&= \frac{1}{\|u\|^p} \Box u;
\end{align*}
since $\nabla_{\R^m} F$ is homogeneous of order $p$, this gives \eqref{nlw} as required.

It remains to establish Theorem \ref{main-3}.
This will be the focus of the remaining sections of the paper.

\section{Eliminating the field}

Having eliminated the potential $F$ from the problem, the next step is (perhaps surprisingly) to eliminate the unknown field $u$, replacing it with quadratic data such as the mass density
\begin{equation}\label{mas}
 M(t,x) := \| u(t,x) \|_{\R^m}^2
\end{equation}
and the energy tensor
\begin{equation}\label{en}
 E_{\alpha \beta}(t,x) := \langle \partial_\alpha u(t,x), \partial_\beta u(t,x) \rangle.
\end{equation}
If $u$ has the discrete self-similarity property \eqref{uts}, then $M$ and $E$ similarly obey the discrete self-similarity properties
\begin{equation}\label{mass-sim}
 M(e^S t, e^S x ) = e^{-\frac{4}{p-1}S} M(t,x)
\end{equation}
and
\begin{equation}\label{dip}
 E_{\alpha \beta}(e^S t, e^S x ) = e^{-\frac{2(p+1)}{p-1}S} E_{\alpha \beta}(t,x).
\end{equation}
Next, observe from the product rule that
\begin{equation}\label{pr}
 \langle u, \Box u \rangle_{\R^m} = \frac{1}{2} \Box M - \eta^{\beta \gamma} E_{\beta \gamma}
\end{equation}
where $\eta$ is the Minkowski metric.  Thus, the defocusing property \eqref{pos} can be rewritten as
\begin{equation}\label{newpos}
\frac{1}{2} \Box M - \eta^{\alpha \beta} E_{\alpha \beta} > 0.
\end{equation}
In a similar spirit, we have
$$ \langle \partial_\alpha u, \Box u \rangle = \partial^\beta E_{\alpha \beta} - \frac{1}{2} \partial_\alpha (\eta^{\beta \gamma} E_{\beta \gamma})$$ 
and hence the equation \eqref{alpha-eq} can be expressed in terms of $M$ and $E$ as
\begin{equation}\label{not}
\partial_\alpha \left( \frac{1}{2} \Box M - \eta^{\beta \gamma} E_{\beta \gamma} \right) = (p+1) \left(\partial^\beta E_{\alpha \beta} - \frac{1}{2} \partial_\alpha (\eta^{\beta \gamma} E_{\beta \gamma})\right).
\end{equation}
Finally, observe that the $2+d \times 2+d$ Gram matrix
\begin{equation}\label{gram}
\begin{pmatrix}
\langle u(t,x), u(t,x) \rangle_{\R^m} & \langle u(t,x), \partial_t u(t,x) \rangle_{\R^m} & \dots & \langle u(t,x), \partial_{x_d} u(t,x) \rangle_{\R^m} \\
\langle \partial_t u(t,x), u(t,x) \rangle_{\R^m} & \langle \partial_t u(t,x), \partial_t u(t,x) \rangle_{\R^m} & \dots & \langle \partial_{t} u(t,x), \partial_{x_d} u(t,x) \rangle_{\R^m} \\
\vdots & \vdots & \ddots & \vdots \\
\langle \partial_{x_d} u(t,x), u(t,x) \rangle_{\R^m} & \langle \partial_{x_d} u(t,x), \partial_{t} u(t,x) \rangle_{\R^m} & \dots & \langle \partial_{x_d} u(t,x), \partial_{x_d} u(t,x) \rangle_{\R^m} 
\end{pmatrix}
\end{equation}
can be expressed in terms of $E,M$ as
\begin{equation}\label{la}
\begin{pmatrix}
M(t,x) & \frac{1}{2} \partial_t M(t,x) & \dots & \frac{1}{2} \partial_{x_d} M(t,x) \\
\frac{1}{2} \partial_t M(t,x) & E_{00}(t,x) & \dots & E_{0d}(t,x) \\
\vdots & \vdots & \ddots & \vdots \\
\frac{1}{2} \partial_{x_d} M(t,x) & E_{d0}(t,x) & \dots & E_{dd}(t,x)
\end{pmatrix}.
\end{equation}
In particular, the matrix \eqref{la} is positive semi-definite for every $t,x$.

It turns out that with the aid of the Nash embedding theorem and our hypothesis that $m$ is large, we can largely reverse the above observations, reducing Theorem \ref{main-3} to the following claim that no longer directly involves the field $u$ (or the range dimension $m$).

\begin{theorem}[Third reduction]\label{main-4}   Let $d = 3$, and let $p > 1 + \frac{4}{d-2}$.   Then there exists $S > 0$ and smooth functions $M\colon \Gamma_d \to \R$ and $E_{\alpha \beta}\colon \Gamma_d \to \R$ for $\alpha,\beta=0,\dots,d$ which are discretely self-similar in the sense of \eqref{mass-sim}, \eqref{dip}, obey the defocusing property \eqref{newpos} and the equation \eqref{not} on $\Gamma_d$ for all $\alpha=0,\dots,d$, and such that the matrix \eqref{la} is strictly positive definite on $\Gamma_d$ (in particular, this forces $M$ to be strictly positive).
\end{theorem}

Let us assume Theorem \ref{main-4} for the moment and show Theorem \ref{main-3}.  Let $d,p,S,M,E_{\alpha \beta}$ be as in Theorem \ref{main-4}, and let $m$ be as in Theorem \ref{main-3}.  Our task is to obtain a function $u\colon \Gamma_d \to \R^m \backslash \{0\}$ obeying all the properties claimed in Theorem \ref{main-3}.

The idea is to build $u$ in such a fashion that \eqref{mas}, \eqref{en} are obeyed.  Accordingly, we will use an ansatz
\begin{equation}\label{ansatz}
 u(t,x) := M(t,x)^{1/2} \theta(t,x)
\end{equation}
for some smooth $\theta\colon \Gamma_d \to S^{m-1}$ to be constructed shortly.  As $M$ is strictly positive, such a function $u$ will be smooth on $\Gamma$ and obey \eqref{mas}; differentiating, we see that
\begin{equation}\label{mas2}
 \langle u, \partial_\alpha u \rangle_{\R^m} = \frac{1}{2} \partial_\alpha M 
\end{equation}
for $\alpha=0,\dots,d$.  If $\theta$ obeys the discrete self-similarity property
\begin{equation}\label{thot}
 \theta( e^S t, e^S x ) = \theta(t,x)
\end{equation}
then $u$ will obey \eqref{uts}.  Thus we shall impose \eqref{thot}, that is to say we assume that $\theta$ is lifted from a smooth map $\tilde \theta\colon \Gamma_d/e^{S\Z} \to S^{m-1}$.

From the product rule and \eqref{mas}, \eqref{mas2} we have (after some calculation)
$$ \langle \partial_\alpha \theta, \partial_\beta \theta \rangle_{\R^m} = M^{-1} \langle \partial_\alpha u, \partial_\beta u \rangle_{\R^m} - M^{-2} (\partial_\alpha M) \partial_\beta M.$$
Thus, if we wish for \eqref{en} to be obeyed, then the $1+d \times 1+d$ Gram matrix
$$ 
(\langle \partial_\alpha \theta, \partial_\beta \theta \rangle_{\R^m} )_{\alpha,\beta=0,\dots,d} $$
must be equal to
\begin{equation}\label{ps}
 M ( E_{\alpha \beta} - (\partial_\alpha M) M^{-1} \partial_\beta M)_{\alpha,\beta=0,\dots,d} .
\end{equation}
The matrix in \eqref{ps} is a Schur complement of the matrix in \eqref{la}.  Since the matrix in \eqref{la} is assumed to be strictly positive definite, we conclude that \eqref{ps} is also.

If we denote the matrix in \eqref{ps} by $g(t,x)$, then  from \eqref{mass-sim}, \eqref{dip} we have the discrete self-similarity property
\begin{equation}\label{gets}
 g(e^S t, e^S x) = e^{-2S} g( t, x ).
\end{equation}
As $g$ is a positive definite and symmetric $1+d \times 1+d$ matrix, we can view $g$ as a smooth Riemannian metric on $\Gamma_d$.
Given that the dilation operator $(t,x) \mapsto (e^S t, e^S x)$ dilates tangent vectors to $\Gamma_d$ by a factor of $e^S$, we see that the metric $g$ is lifted from a smooth Riemannian metric $\tilde g$ on the quotient space $\Gamma_d / e^{S\Z}$.

The space $(\Gamma_d / e^{S\Z},\tilde g)$ is a smooth compact $1+d$-dimensional Riemannian manifold with boundary; it is easy to embed it in a smooth compact $1+d$-dimensional Riemannian manifold without boundary (for instance by using the theorems in \cite{seeley}).  Applying the Nash embedding theorem (for instance in the form in \cite{gunther}), we can thus isometrically embed $(\Gamma_d / e^{S\Z},\tilde g)$ in a Euclidean space $\R^D$ with $D := \max(\frac{(d+1)(d+6)}{2}, \frac{(d+1)(d+4)}{2}+5)$.  The embedded copy of $(\Gamma_d / e^{S\Z},\tilde g)$ is compact and is thus contained in a cube $[-R,R]^D$ for some finite $R$.  We use a generic\footnote{We thank Marc Nardmann for this argument, which improved the value of $m$ from our previous argument by a factor of approximately two.} linear isometry from $\R^D$ to $\R^{D+1}$ to embed $[-R,R]^D$ to some compact subset of $\R^{D+1}$.  The image of this isometry is a generic hyperplane, which can be chosen to avoid the lattice $\frac{1}{\sqrt{2D+2}} \Z^{D+1}$, and thus we can embed $[-R,R]^D$ isometrically into the torus $\R^{D+1} / \frac{1}{\sqrt{D+1}} \Z^{D+1}$, which is isometric to $\frac{1}{\sqrt{D+1}} (S^1)^{D+1}$.  But from Pythagoras' theorem, $\frac{1}{\sqrt{D+1}} (S^1)^{D+1}$ is contained in $S^{2D+1}$, which is in turn contained in $S^{m-1}$ by the largeness hypothesis on $m$.  Thus we have an isometric embedding $\tilde \theta\colon \Gamma_d / e^{S\Z} \to S^{m-1}$ from $(\Gamma_d/e^{S\Z}, \tilde g)$ into the round sphere $S^{m-1}$.  In particular, $\tilde \theta$ is injective and immersed, and lifting $\tilde \theta$ back to $\Gamma_d$, we obtain a smooth map $\theta\colon \Gamma_d \to S^{m-1}$ with Gram matrix \eqref{ps} that is discretely self-similar in the sense of \eqref{thot}, so that the function $u$ defined by \eqref{ansatz} obeys \eqref{uts}.  Reversing the calculations that led to \eqref{ps}, we see that the Gram matrix \eqref{gram} of $u$ is given by \eqref{la}.  In particular, \eqref{en} holds.  Reversing the derivation of \eqref{newpos}, we now obtain \eqref{pos}, while from reversing the derivation of \eqref{not}, we obtain \eqref{alpha-eq}.  We have now obtained all the required properties claimed by Theorem \ref{main-3}, as desired.

It remains to establish Theorem \ref{main-4}.
This will be the focus of the remaining sections of the paper.

\section{Reduction to a self-similar $1+1$-dimensional problem}

In reducing Theorem \ref{main} to Theorem \ref{main-4}, we have achieved the somewhat remarkable feat of converting a nonlinear PDE problem to a convex (or positive semi-definite) PDE problem, in that all of the constraints\footnote{Compare with the ``kernel trick'' in machine learning, or with semidefinite relaxation in optimization.} on the remaining unknowns $M, E_{\alpha \beta}$ are linear equalities and inequalities, or assertions that certain matrices are positive definite.  Among other things, this shows that if one has a given solution $M, E_{\alpha \beta}$ to Theorem \ref{main-4}, and then one averages that solution over some compact symmetry group that acts on the space of such solutions, then the average will also be a solution to Theorem \ref{main-4}.  In particular, one can then reduce without any loss of generality to considering solutions that are invariant with respect to that symmetry.

For instance, given that $M, E_{\alpha \beta}$ are already discretely self-similar by \eqref{mass-sim}, \eqref{dip}, the space of solutions has an action of the compact dilation group $\R^+ / e^{S\Z}$, with (the quotient representative of) any real number $\lambda > 0$ acting on $M, E_{\alpha \beta}$ by the action
$$ (\lambda \cdot M)(t,x) := \frac{1}{\lambda^{\frac{4}{p-1}}} M( \frac{t}{\lambda}, \frac{x}{\lambda})$$
and
$$ (\lambda \cdot E_{\alpha \beta})(t,x) := \frac{1}{\lambda^{\frac{2(p+1)}{p-1}}} E_{\alpha \beta}( \frac{t}{\lambda}, \frac{x}{\lambda});$$
this is initially an action of the multiplicative group $\R^+$, but descends to an action of $\R^+ / e^{S\Z}$ thanks to \eqref{mass-sim}, \eqref{dip}.  By the preceding discussion, we may restrict without loss of generality to the case when $M, E_{\alpha\beta}$ are invariant with respect to this $\R^+ / e^{S\Z}$, or equivalently that $M$ and $E_{\alpha \beta}$ are homogeneous of order $-\frac{4}{p-1}$ and $\frac{2(p+1)}{p-1}$ respectively.  With this restriction, the parameter $S$ no longer plays a role and may be discarded.

\begin{remark} This reduction may seem at first glance to be in conflict with the negative result in Proposition \ref{nose}.  However, the requirement that the mass density $M$ and the energy tensor $E_{\alpha \beta}$ be homogeneous is strictly weaker than the hypothesis that the field $u$ itself is homogeneous.  For instance, one could imagine a ``twisted self-similar'' solution in which the homogeneity condition \eqref{homog} on $u$ is replaced with a more general condition of the form
$$ u( \lambda t, \lambda x ) = \lambda^{-\frac{2}{p-1}} \exp( J \log \lambda ) u(t,x) $$
for all $(t,x) \in \Gamma_d$ and $\lambda > 0$, where $J\colon \R^m \to \R^m$ is a fixed skew-adjoint linear transformation.  (To be compatible with \eqref{nlw}, one would also wish to require that the potential $F$ is invariant with respect to the orthogonal transformations $\exp( s J)$ for $s \in \R$.)  Such solutions $u$ would not be homogeneous, but the associated densities $M, E_{\alpha \beta}$ would still be homogeneous of the order specified above.  
\end{remark}

We may similarly apply the above reductions to the orthogonal group $O(d)$, which acts on the scalar field $M$ and on the $2$-tensor $E_{\alpha \beta}$ in the usual fashion, thus
$$ (U M)(t,x) := M( t, U^{-1} x )$$
and
$$ (U E)_{\alpha \beta}(t,x) (Uv)^\alpha (Uv)^\beta = E_{\alpha \beta}(t,U^{-1} x) v^\alpha v^\beta$$
for all $(t,x) \in \Gamma_d$, $U \in O(d)$, and $v \in \R^{1+d}$, where $U$ acts on $\R^{1+d}$ by $(t,x) \mapsto (t,Ux)$.  This allows os to reduce to fields $M, E_{\alpha \beta}$ which are $O(d)$-invariant, thus $M$ is spherically symmetric, and $E_{\alpha \beta}$ takes the form\footnote{To see that $E_{\alpha \beta}$ must be of this form, rotate the spatial variable $x$ to equal $x = r e_1$, then use the orthogonal transformation $(x_1,x_2,\dots,x_d) \mapsto (x_1,-x_2,\dots,-x_d)$, which preserves $re_1$, to see that $E_{0i} = E_{1i} = 0$ for all $i=2,\dots,d$; further use of orthogonal transformations preserving $re_1$ can be then used to show that $E_{ij}=0$ and $E_{ii}=E_{jj}$ for $2 \leq i < j \leq d$ (basically because the only matrices that commute with all orthogonal transformations are scalar multiples of the identity).  This places $E_{\alpha \beta}$ in the desired form in the $x=re_1$ case, and the general case follows from rotation.}
\begin{align}
E_{00} &= E_{tt} \label{et1} \\
E_{0i} = E_{i0} &= \frac{x_i}{r} E_{tr}\label{et2} \\
E_{ij} &= \frac{x_i x_j}{r^2} (E_{rr} - E_{\omega \omega}) + \delta_{ij} E_{\omega \omega} \label{et3}
\end{align}
for $i,j = 1,\dots,d$ and some spherically symmetric scalar functions $E_{tt}, E_{tr}, E_{rr}, E_{\omega \omega}$, where $r := |x|$ is the radial variable and $\delta_{ij}$ is the Kronecker delta.  Observe that if $E_{tt}, E_{rr}, E_{\omega\omega}:\Gamma_1 \to \R$ are smooth even functions and $E_{tr}\colon \Gamma_1 \to \R$ a smooth odd function on the $1+1$-dimensional light cone
$$ \Gamma_1 := \{ (t,r) \in \R^{1+1}: t > 0; -t \leq r \leq t \}$$
with $E_{rr}-E_{\omega\omega}$ vanishing to second order at $r=0$, then the above equations define a smooth field $E_{\alpha \beta}$ on $\Gamma_d$, which will be homogeneous of order $-\frac{2(p+1)}{p-1}$ if $E_{tt}, E_{tr}, E_{rr}, E_{\omega\omega}$ are.

Using polar coordinates, we have
$$
\frac{1}{2} \Box M - \eta^{\beta \gamma} E_{\beta \gamma} 
= \frac{1}{2} \left(-\partial_{tt} M + \partial_{rr} M + \frac{d-1}{r} M\right) - (-E_{tt} + E_{rr} + (d-1) E_{\omega\omega}) $$
thus the condition \eqref{newpos} is now
\begin{equation}\label{newerpos}
 \frac{1}{2} \left(-\partial_{tt} M + \partial_{rr} M + \frac{d-1}{r} M\right) - (-E_{tt} + E_{rr} + (d-1) E_{\omega\omega}) > 0.
\end{equation}
By rotating $x$ to be of the form $x=re_1$, we see that the matrix \eqref{la} is conjugate to
$$
\begin{pmatrix}
M & \frac{1}{2} \partial_t M & \frac{1}{2} \partial_r M & 0 & \dots & 0 \\
\frac{1}{2} \partial_t M & E_{tt} & E_{tr} & 0 & \dots & 0 \\
\frac{1}{2} \partial_r M & E_{tr} & E_{rr} & 0 & \dots & 0 \\
0 & 0 & 0 & E_{\omega\omega} & \dots & 0 \\
\vdots & \vdots & \vdots & \vdots & \ddots & \vdots \\
0 & 0 & 0 & 0 & \dots & E_{\omega\omega}
\end{pmatrix}
$$
so the positive-definiteness of \eqref{la} is equivalent to the positive definiteness of the $3 \times 3$ matrix
\begin{equation}\label{la-2}
\begin{pmatrix}
M & \frac{1}{2} \partial_t M & \frac{1}{2} \partial_r M \\
\frac{1}{2} \partial_t M & E_{tt} & E_{tr} \\
\frac{1}{2} \partial_r M & E_{tr} & E_{rr} 
\end{pmatrix}
\end{equation}
together with the positivity of $E_{\omega\omega}$.  It will be convenient to isolate the $r=0$ case of this condition (in order to degenerate $E_{rr}$ to zero at $r=0$ later in the argument).  In this case, the odd functions $\partial_r M$ and $E_{tr}$ vanish, and $E_{rr}$ is equal to $E_{\omega\omega}$, so the condition reduces to the positive definiteness of the $2 \times 2$ matrix
\begin{equation}\label{la-3}
\begin{pmatrix}
M & \frac{1}{2} \partial_t M \\
\frac{1}{2} \partial_t M & E_{tt} 
\end{pmatrix}
\end{equation}
together with the aforementioned positivity of $E_{\omega\omega}$.

Finally, we turn to the condition \eqref{not}.  Again, we can rotate the position $x$ to be of the form $x=re_1$.  In the angular cases $\alpha=2,\dots,d$, both sides of \eqref{not} automatically vanish, basically because $\partial_\alpha f(re_1)=0$ for any spherically symmetric $f$ (and because $E_{\alpha \beta}$ vanishes to second order for any $\beta \neq \alpha$).  So the only non-trivial cases of \eqref{not} are $\alpha=0$ and $\alpha=1$.  Applying \eqref{et1}, \eqref{et2}, \eqref{et3}, we can write these cases of \eqref{not} as
\begin{equation}\label{not-0}
\begin{split}
&\partial_t \left[\frac{1}{2} \left(-\partial_{tt} M + \partial_{rr} M + \frac{d-1}{r} M\right) - (-E_{tt} + E_{rr} + (d-1) E_{\omega\omega})\right]\\
&\quad = (p+1) \left[-\partial_t E_{tt} + \partial_r E_{tr} + \frac{d-1}{r} E_{tr} - \frac{1}{2} \partial_t(-E_{tt} + E_{rr} + (d-1) E_{\omega\omega}) \right]
\end{split}
\end{equation}
and
\begin{equation}\label{not-1}
\begin{split}
&\partial_r \left[\frac{1}{2} \left(-\partial_{tt} M + \partial_{rr} M + \frac{d-1}{r} M\right) - (-E_{tt} + E_{rr} + (d-1) E_{\omega\omega})\right]\\
&\quad = (p+1) \left[-\partial_t E_{tr} + \partial_r E_{rr} + \frac{d-1}{r} (E_{rr}-E_{\omega\omega}) - \frac{1}{2} \partial_r(-E_{tt} + E_{rr} + (d-1) E_{\omega\omega}) \right]
\end{split}
\end{equation}
respectively.

To summarise, we have reduced Theorem \ref{main-4} to

\begin{theorem}[Fourth reduction]\label{main-5}   Let $d = 3$, and let $p > 1 + \frac{4}{d-2}$.   Then there exist smooth even functions $M, E_{tt}, E_{rr}, E_{\omega\omega}\colon \Gamma_1 \to \R$ and a smooth odd function $E_{tr}\colon \Gamma_1 \to \R$, with $M$ homogeneous of order $-\frac{4}{p-1}$ and $E_{tt}, E_{tr}, E_{rr}, E_{\omega\omega}$ homogeneous of order $-\frac{2(p+1)}{p-1}$, and with $E_{rr}-E_{\omega\omega}$ vanishes to second order at $r=0$, obeying the defocusing property \eqref{newerpos} and the equations \eqref{not-0}, \eqref{not-1} on $\Gamma_1$, such that
\begin{equation}\label{eww}
E_{\omega\omega} > 0
\end{equation}
and the $3 \times 3$ matrix \eqref{la-2} is strictly positive definite on $\Gamma_1$ with $r \neq 0$, and the $2 \times 2$ matrix \eqref{la-3} is positive definite when $r=0$.
\end{theorem}

It remains to prove Theorem \ref{main-5}.  To do so, we make a few technical relaxations.  Firstly, we claim that we may relax the strict conditions \eqref{newerpos} and \eqref{eww} to their non-strict counterparts
\begin{equation}\label{newerpos-weak}
 \frac{1}{2} \left(-\partial_{tt} M + \partial_{rr} M + \frac{d-1}{r} M\right) - (-E_{tt} + E_{rr} + (d-1) E_{\omega\omega}) \geq 0
\end{equation}
and
\begin{equation}\label{eww-weak}
E_{\omega\omega} \geq 0.
\end{equation}
To see this, suppose that $M, E_{tt}, E_{tr}, E_{rr}, E_{\omega\omega}$ obey the conclusions of Theorem \ref{main-5} with the conditions \eqref{newerpos}, \eqref{eww} replaced by \eqref{newerpos-weak}, \eqref{eww-weak}.  We let $\eps>0$ be a small quantity to be chosen later, and define new fields $M^\eps, E^\eps_{tt}, E^\eps_{tr}, E^\eps_{rr}, E^\eps_{\omega\omega}$ by the formulae\footnote{The ability to freely manipuate the fields $M, E_{tt}, E_{tr}, E_{rr}, E_{\omega\omega}$ in this fashion is a major advantage of the formulation of Theorem \ref{main-5}.  It would be very difficult to perform analogous manipulations if the original field $u$ or the potential $F$ were still present.}
\begin{align*}
M^\eps &:= M - c \eps t^{-\frac{4}{p-1}} \\
E^\eps_{tt} &:= E_{tt} - (d+1) \eps t^{-\frac{2(p+1)}{p-1}} \\
E^\eps_{tr} &:= E_{tr} \\
E^\eps_{rr} &:= E_{rr} + \eps t^{-\frac{2(p+1)}{p-1}} \\
E^\eps_{\omega\omega} &:= E_{\omega\omega} + \eps t^{-\frac{2(p+1)}{p-1}}
\end{align*}
where $c$ is the constant such that
$$ \frac{1}{2} c \frac{4}{p-1} \frac{p+3}{p-1} - (2d+1) = \frac{p+1}{2}.$$
Clearly these new fields $M^\eps, E^\eps_{tt}, E^\eps_{tr}, E^\eps_{rr}, E^\eps_{\omega\omega}$ are still smooth, with $M^\eps, E^{\eps}_{tt}, E^{\eps}_{rr}, E^\eps_{\omega\omega}$ even and $E^\eps_{tr}$ odd, with $M^\eps$ homogeneous of order $-\frac{4}{p-1}$ and $E^\eps_{tt}, E^\eps_{tr}, E^\eps_{rr}, E^\eps_{\omega\omega}$ homogeneous of order $-\frac{2(p+1)}{p-1}$, with $E^\eps_{rr}-E^\eps_{\omega\omega}$ vanishing to second order at $r=0$.  A calculation using the definition of $c$ shows that the equations \eqref{not-0}, \eqref{not-1} continue to be obeyed when the fields $M, E_{tt}, E_{tr}, E_{rr}, E_{\omega\omega}$ are replaced by $M^\eps, E^\eps_{tt}, E^\eps_{tr}, E^\eps_{rr}, E^\eps_{\omega\omega}$.  With this replacement, the left-hand side of \eqref{newerpos-weak} increases by $\frac{p+1}{2} \eps t^{-\frac{2(p+1)}{p-1}}$, and so \eqref{newerpos} now holds.  The remaining task is to show that with these new fields $M^\eps, E^\eps_{tt}, E^\eps_{tr}, E^\eps_{rr}, E^\eps_{\omega\omega}$, \eqref{la-2} is positive definite when $r \neq 0$ and \eqref{la-3} is positive definite when $r=0$.  By the scale invariance it suffices to verify these latter properties when $t=1$.  The positive definiteness of \eqref{la-3} when $r=0$ then follows by continuity for $\eps$ small enough.  For \eqref{la-2}, we have to take a little care because the condition $r \neq 0$ is non-compact.  We need to ensure the positive definiteness of 
$$
\begin{pmatrix}
M-c\eps & \frac{1}{2} \partial_t M + \frac{2c}{p-1} \eps & \frac{1}{2} \partial_r M \\
\frac{1}{2} \partial_t M + \frac{2c}{p-1} \eps & E_{tt} - (d+1)\eps & E_{tr} \\
\frac{1}{2} \partial_r M & E_{tr} & E_{rr}+\eps 
\end{pmatrix}
$$
when $t=1$ and $r \neq 0$, for $\eps$ small enough.  Continuity will ensure this if $|r|$ is bounded away from zero (independently of $\eps$), so we may assume that $r$ is in a small neighbourhood of the origin (independent of $\eps$).  Given that the above matrix is already positive definite when $\eps = 0$, it suffices by a continuity argument to show that the above matrix has positive determinant for sufficiently small $\eps$; by the hypothesis \eqref{la-2} and the fundamental theorem of calculus, it thus suffices to show that
$$
\frac{d}{d\eps} \operatorname{det}
\begin{pmatrix}
M-c\eps & \frac{1}{2} \partial_t M + \frac{2c}{p-1} \eps & \frac{1}{2} \partial_r M \\
\frac{1}{2} \partial_t M + \frac{2c}{p-1} \eps & E_{tt} - (d+1)\eps & E_{tr} \\
\frac{1}{2} \partial_r M & E_{tr} & E_{rr}+\eps 
\end{pmatrix} > 0 $$
for $r$ near zero and sufficiently small $\eps$.  But since $\partial_r M, E_{tr}, E_{rr}$ vanish at $r=0$, we can use cofactor expansion to write the left-hand side as
$$
\operatorname{det} \begin{pmatrix} M(1,0) & \frac{1}{2} \partial_t M(1,0) \\ \frac{1}{2} \partial_t M(1,0) & E_{tt}(1,0) \end{pmatrix} + O( |r| ) + O(\eps) 
$$
and the claim then follows from the hypothesis \eqref{la-2}.  This concludes the relaxation of the conditions \eqref{newerpos}, \eqref{eww} to \eqref{newerpos-weak}, \eqref{eww-weak}.

Now that we allow equality in \eqref{eww-weak}, we sacrifice some generality by restricting to the special case $E_{\omega\omega} = 0$ (which basically corresponds to considering spherically symmetric blowup solutions).  While this gives up some flexibility, this will simplify our calculations a bit as we now only have four fields $M, E_{tt}, E_{tr}, E_{rr}$ to deal with, rather than five.

Until now we have avoided using the hypothesis $d=3$.  Now we will embrace this hypothesis.  In Proposition \ref{nose} it was convenient to make the change of variables $\phi=r^{(d-1)/2} u = ru$ to eliminate lower order terms such as $\frac{d-1}{r} \partial_r u$; this change of variables is particularly pleasant in the three-dimensional case as the lower-order term involving the coefficient $\frac{(d-1)(d-3)}{4}$ vanishes completely (this vanishing is closely tied to the strong Huygens principle in three dimensions).  The corresponding change of variables in this setting, aimed at eliminating the lower order terms $\frac{d-1}{r} E_{tr}, \frac{d-1}{r} E_{rr}$ in \eqref{not-0}, \eqref{not-1}, is to replace the fields $M, E_{tt}, E_{tr}, E_{rr}$ by the fields $\tilde M, \tilde E_{tt}, \tilde E_{tr}, \tilde E_{rr}\colon \Gamma_1 \to \R^+$ defined by
\begin{align*}
\tilde M &:= r^2 M \\
\tilde E_{tt} &:= r^2 E_{tt} \\
\tilde E_{tr} &:= r^2 E_{tr} + \frac{1}{2} r \partial_t M \\
&= r^2 E_{tr} + \frac{1}{2r} \partial_t \tilde M \\
\tilde E_{rr} &:= r^2 E_{rr} + r \partial_r M + M\\
&= r^2 E_{rr} + \frac{1}{r} \partial_r \tilde M - \frac{1}{r^2} \tilde M.
\end{align*}
Observe that if $\tilde M, \tilde E_{tt}, \tilde E_{rr}$ are smooth and even, and $\tilde E_{tr}$ odd, with $\tilde M, \tilde E_{tt}$ vanishing to second order at $r=0$, $\tilde E_{tr} - \frac{1}{2r} \partial_t \tilde M$ vanishing to third order, and $\tilde E_{rr} - \frac{1}{r} \partial_r \tilde M + \frac{1}{r^2} \tilde M$ to fourth order, then these fields determine smooth fields $M, E_{tt}, E_{tr}, E_{rr}$ with $M, E_{tt}, E_{rr}$ even, $E_{tr}$ odd, and $E_{rr}$ vanishing to second order at $r=0$.  Furthermore, if $\tilde M$ is homogeneous of order $\frac{2p-6}{p-1}$ and $\tilde E_{tt}, \tilde E_{tr}, \tilde E_{rr}$ are homogeneous of order $-\frac{4}{p-1}$, then $M$ will be homogeneous of order $-\frac{4}{p-1}$ and $E_{tt}, E_{tr}, E_{rr}$ will be homogeneous of order $-\frac{2(p+1)}{p-1}$.

If we introduce the quantity
\begin{equation}\label{vdef}
 V := \frac{1}{p+1} \left( \frac{1}{2} ( - \partial_{tt} \tilde M + \partial_{rr} \tilde M ) + \tilde E_{tt} - \tilde E_{rr} \right) 
\end{equation}
then a brief calculation shows that
$$ V = \frac{r^2}{2 (p+1)} \left((-\partial_{tt} M + \partial_{rr} M + \frac{2}{r} M) - (-E_{tt} + E_{rr})\right) $$
and so the condition \eqref{newerpos-weak} is equivalent to
\begin{equation}\label{rpos}
 V \geq 0. 
\end{equation}
The equations \eqref{not-0}, \eqref{not-1} can now be expressed as
$$
\partial_t \left[\frac{1}{r^2} V\right]
= -\partial_t E_{tt} + \partial_r E_{tr} + \frac{2}{r} E_{tr} - \frac{1}{2} \partial_t(-E_{tt} + E_{rr}) 
$$
and
$$
\partial_r \left[\frac{1}{r^2} V\right]
= -\partial_t E_{tr} + \partial_r E_{rr} + \frac{2}{r} E_{rr} - \frac{1}{2} \partial_r(-E_{tt} + E_{rr})
$$
which rearrange as an energy conservation law
$$ \partial_t \left( \frac{1}{2} E_{tt} + \frac{1}{2} E_{rr} + \frac{1}{r^2} V \right) = \partial_r E_{tr} + \frac{2}{r} E_{tr}$$
and a momentum conservation law
$$ \partial_t E_{tr} = \partial_r \left( \frac{1}{2} E_{tt} + \frac{1}{2} E_{rr} - \frac{1}{r^2} V \right) + \frac{2}{r} E_{rr};$$
multiplying these equations by $r^2$ and writing $E_{tt}, E_{tr}, E_{rr}$ in terms of $\tilde E_{tt}, \tilde E_{tr}, \tilde E_{rr}$ and $\tilde M$ one obtains (after some calculation, as well as \eqref{vdef} in the case of \eqref{notnew-1}) the slightly simpler equations
\begin{equation}\label{notnew-0}
\partial_t \left( \frac{1}{2} \tilde E_{tt} + \frac{1}{2} \tilde E_{rr} + V \right) = \partial_r \tilde E_{tr} 
\end{equation}
and
\begin{equation}\label{notnew-1}
\partial_t \tilde E_{tr} = \partial_r ( \frac{1}{2} \tilde E_{tt} + \frac{1}{2} \tilde E_{rr} - V ) - \frac{p-1}{r} V.
\end{equation}
The expressions in \eqref{notnew-0} are even, while the expressions in \eqref{notnew-1} are odd.  Thus we may combine these equations into a single equation by adding them together, which after some rearranging becomes the transport type equation
\begin{equation}\label{storp}
 (\partial_t - \partial_r) e_+ + (\partial_t + \partial_r) V = - \frac{p-1}{r} V 
\end{equation}
where $e_+$ is the null energy density
\begin{equation}\label{null-energy}
e_+ := \frac{1}{2} \tilde E_{tt} + \frac{1}{2} \tilde E_{rr} + \tilde E_{tr}.
\end{equation}

\begin{remark} It may be instructive to derive these equations in the specific context of a solution $u\colon \Gamma_3 \to \R$ to the scalar defocusing NLW
$$ \Box u = |u|^{p-1} u $$
which in polar coordinates becomes
$$ -\partial_{tt} u + \partial_{rr} u + \frac{2}{r} u = |u|^{p-1} u.$$
Making the change of variables $\phi = ru$, this becomes
$$ -\partial_{tt} \phi + \partial_{rr} \phi = \frac{|\phi|^{p-1} \phi}{r^{p-1}}.$$
Introducing the null energy
$$ e_+ := \frac{1}{2} |\partial_t \phi + \partial_r \phi|^2$$
and the potential energy
$$ V := \frac{1}{p+1} \frac{|\phi|^{p+1}}{r^{p-1}}$$
as well as the additional densities
$$ \tilde M := |\phi|^2; \quad \tilde E_{tt} := |\partial_t \phi|^2; \quad \tilde E_{rr} := |\partial_r \phi|^2; \quad \tilde E_{tr} := \partial_t \phi \partial_r \phi,$$
one can readily verify the identities \eqref{vdef}, \eqref{storp}, and \eqref{null-energy}.  Similarly for the other properties of $\tilde M, \tilde E_{tt}, \tilde E_{rr}, \tilde E_{tr}$ identified in this section.
\end{remark}

Finally, we translate the positive definiteness of \eqref{la-2} (when $r \neq 0$) and \eqref{la-3} (when $r=0$) into conditions involving the fields $\tilde M, \tilde E_{tt}, \tilde E_{rr}, \tilde E_{tr}$.  From the identity
$$
\begin{pmatrix}
\tilde M & \frac{1}{2} \partial_t \tilde M & \frac{1}{2} \partial_r \tilde M \\
\frac{1}{2} \partial_t \tilde M & \tilde E_{tt} & \tilde E_{tr} \\
\frac{1}{2} \partial_r \tilde M & \tilde E_{tr} & \tilde E_{rr}
\end{pmatrix}
= r^2
\begin{pmatrix}
1 & 0 & 0 \\
0 & 1 & 0 \\
1/r & 0 & 1
\end{pmatrix}
\begin{pmatrix}
M & \frac{1}{2} \partial_t M & \frac{1}{2} \partial_r M \\
\frac{1}{2} \partial_t M & E_{tt} & E_{tr} \\
\frac{1}{2} \partial_r M & E_{tr} & E_{rr}
\end{pmatrix}
\begin{pmatrix}
1 & 0 & 1/r \\
0 & 1 & 0 \\
0 & 0 & 1
\end{pmatrix}
$$
we see (for $r \neq 0$) that \eqref{la-2} is strictly positive-definite if and only if the matrix
\begin{equation}\label{la-4}
\begin{pmatrix}
\tilde M & \frac{1}{2} \partial_t \tilde M & \frac{1}{2} \partial_r \tilde M \\
\frac{1}{2} \partial_t \tilde M & \tilde E_{tt} & \tilde E_{tr} \\
\frac{1}{2} \partial_r \tilde M & \tilde E_{tr} & \tilde E_{rr}
\end{pmatrix}
\end{equation}
is strictly positive-definite.  Now we turn to \eqref{la-3} when $r=0$.  By homogeneity, it suffices to verify this condition when $(t,r) = (1,0)$.  From \eqref{homog}, we have $\partial_t M(1,0) = - \frac{4}{p-1} M(1,0)$, so the positive definiteness of \eqref{la-3} is equivalent to the condition
$$ E_{tt}(1,0) > \left(\frac{2}{p-1}\right)^2 M(1,0) > 0 $$
which in terms of $\tilde E_{tt}, \tilde M$ becomes
\begin{equation}\label{mur}
\partial_{rr} \tilde E_{tt}(1,0) > \left(\frac{2}{p-1}\right)^2 \partial_{rr} \tilde M(1,0) > 0.
\end{equation}
Summarising the above discussion, we now see that Theorem \ref{main-5} is a consequence of the following statement.

\begin{theorem}[Fifth reduction]\label{main-6}   Let $p > 5$.  Then there exist smooth even functions $\tilde M, \tilde E_{tt}, \tilde E_{rr}\colon \Gamma_1 \to \R$ and a smooth odd function $\tilde E_{tr}\colon \Gamma_1 \to \R$, with $\tilde M$ homogeneous of order $\frac{2p-6}{p-1}$ and $\tilde E_{tt}, \tilde E_{tr}, \tilde E_{rr}$ homogeneous of order $-\frac{4}{p-1}$, with $\tilde M, \tilde E_{tt}$ vanishing to second order at $r=0$, $\tilde E_{tr} - \frac{1}{2r} \partial_t \tilde M$ vanishing to third order, and $\tilde E_{rr} - \frac{1}{r} \partial_r \tilde M + \frac{1}{r^2} \tilde M$ to fourth order.  Furthermore, if one defines the fields $V,e_+\colon \Gamma_1 \to \R$ by \eqref{vdef} and \eqref{null-energy}, we have the weak defocusing property \eqref{rpos} and the null transport equation \eqref{storp}.  Finally, the matrix \eqref{la-4} is strictly positive definite for $r \neq 0$, and for $r=0$ one has the condition \eqref{mur}.
\end{theorem}

It remains to establish Theorem \ref{main-6}. This will be the focus of the final section of the paper.

\section{Constructing the mass and energy fields}

Fix $p>5$.  We will need a large constant $A>1$ depending only on $p$, and then sufficiently small parameter $\delta>0$ (depending on $p,A$) to be chosen later.  We use the notation $X \lesssim Y$, $Y \gtrsim X$, or $X = O(Y)$ to denote an estimate of the form $|X| \leq CY$, where $C$ can depend on $p$ but is independent of $\delta, A$.  

We need to construct smooth fields $\tilde M, \tilde E_{tt}, \tilde E_{rr}, \tilde E_{tr}\colon \Gamma_1 \to \R$ which generate some further fields $V, e_+\colon \Gamma_1 \to \R$, which are all required to obey a certain number of constraints.  The problem is rather underdetermined, and so there will be some flexibility in selecting these fields; most of these fields will end up being concentrated in the region $\{ (t,r) \in \Gamma_1: r = (\pm 1 + O(\delta)) t \}$ near the boundary of the light cone.  Given that the constraint \eqref{storp} only involves the two fields $V$ and $e_+$, it is natural to proceed by constructing $V$ and $e_+$ first.  In fact we will proceed as follows.

\subsection{Selection of $e_+$ in the left half of the cone}   

We begin by making a choice for the function $e_+\colon \Gamma_1 \to \R$ in the left half $\Gamma_1^l := \{ (t,r) \in \Gamma_1: r \leq 0 \}$ of the cone.  When $t=1$, we choose $e_+(1,r)$ to be a smooth function with the following properties:
\begin{itemize}
\item One has
\begin{equation}\label{ea}
 e_+(1,r) = (1+r)^{-\frac{4}{p-1}}
\end{equation}
for $-1+\delta \leq r \leq 0$.
\item One has
\begin{equation}\label{eb}
e_+(1,r) \geq (1+r)^{-\frac{4}{p-1}}
\end{equation}
for $-1+\frac{\delta}{2} \leq r \leq -1+\delta$.  Furthermore, one has
\begin{equation}\label{big}
\int_{-1+\frac{\delta}{2}}^{-1+\delta} e_+(1,r)\ dr \geq A \delta^{1-\frac{4}{p-1}}.
\end{equation}
\item One has
\begin{equation}\label{ec}
\delta^{-\frac{4}{p-1}} \lesssim e_+(1,r) \lesssim A \delta^{-\frac{4}{p-1}} 
\end{equation}
and
\begin{equation}\label{ed}
|\frac{d}{dr} e_+(1,r)| \lesssim A \delta^{-\frac{p+3}{p-1}} 
\end{equation}
for $-1 \leq r \leq -1+\delta$.
\end{itemize}
Clearly we can find a smooth function $r \to e_+(1,r)$ on $[-1,0]$ with these properties.  We then extend $e_+$ to the entire left half $\Gamma_1^l$ of the cone by requiring it to be homogeneous of order $-\frac{4}{p-1}$, thus
\begin{equation}\label{etr-def}
e_+(t,r) := t^{-\frac{4}{p-1}} e_+(1, \frac{r}{t}).
\end{equation}
In particular, $e_+$ is smooth on this half of the cone, and we have
$$ e_+(t,r) = (t+r)^{-\frac{4}{p-1}}$$
for $-(1-\delta)t \leq r \leq 0$.

The properties \eqref{ea}-\eqref{ed} are largely used to ensure that the potential energy $V$ that we will construct below is non-negative.

\subsection{Selection of $V$ in the left half of the cone}  

Once $e_+$ has been selected on $\Gamma_1^l$, we construct $V$ on $\Gamma_1^l$ by solving \eqref{storp}, or more explicitly by the formula
\begin{equation}\label{vater}
V(t,r) := \frac{1}{2|r|^{p-1}} \int_r^0 |s|^{p-1}  ((\partial_t - \partial_r) e_+) (t-r+s, s)\ ds
\end{equation}
for $-t \leq r < 0$.  Note that as $(\partial_t - \partial_r) e_+$ vanishes for $-(1-\delta)t < r < 0$, $V$ vanishes on this region also, and so one can smoothly extend $V$ to all of $\Gamma_1^l$.  It is easy to see that $V$ is homogeneous of order $-\frac{4}{p-1}$.  From the fundamental theorem of calculus and the chain rule, we have
$$
(\partial_t + \partial_r) (|r|^{p-1} V) = |r|^{p-1} (\partial_t - \partial_r) e_+ $$
for $-t \leq r < 0$, and hence by the product rule we see that \eqref{storp} is obeyed for $-t \leq r < 0$, and hence to all of $\Gamma^1_l$ by smoothness.

We have already seen that $V$ vanishes in the region $-(1-\delta)t < r \leq 0$.  In the region $-t \leq r \leq -(1-\delta)t$, we have the following estimate and non-negativity property:

\begin{proposition}\label{vpos}  For $-t \leq r \leq -(1-\delta)t$, we have
$$ 0 \leq V(t,r) \lesssim A t^{-\frac{4}{p-1}} \delta^{\frac{p-5}{p-1}}.$$
\end{proposition}

We remark that to get the lower bound $V(t,r)$, the supercriticality hypothesis $p>5$ will be crucial.

\begin{proof}  By homogeneity we may assume that $t-r=2$, so that $t = 1 - O(\delta)$, $r = -1 + O(\delta)$, and it will suffice to show that
\begin{equation}\label{dipp}
 0 \leq V(t,r) \lesssim A \delta^{\frac{p-5}{p-1}}.
\end{equation}
Write $e_+(t,r) = (t+r)^{-\frac{4}{p-1}} + f(t,r)$, then from \eqref{vater} we have
\begin{equation}\label{vater-2}
V(t,r) = \frac{1}{2|r|^{p-1}} \int_r^0 |s|^{p-1} ((\partial_t - \partial_r) f) (2+s, s)\ ds.
\end{equation}
The function $f$ is homogeneous of order $-\frac{4}{p-1}$, hence by \eqref{euler}
$$
(t \partial_t + t \partial_r ) f = -\frac{4}{p-1} f.$$
From the identity
$$
\partial_t - \partial_r = - \frac{t+r}{t-r} (\partial_t + \partial_r) + \frac{2}{t-r} (t \partial_t + t \partial_r )$$
and the chain rule, we thus have
$$
((\partial_t - \partial_r) f) (2+s, s) = - (1+s) \frac{d}{ds} f(2+s,s) - \frac{4}{p-1} f(2+s,s).$$
Inserting this into \eqref{vater-2} and integrating by parts, we conclude that
$$ V(t,r) = \frac{1+r}{2} f(t,r) + \frac{1}{2|r|^{p-1}} \int_r^0 \frac{d}{ds}( |s|^{p-1} (1+s) ) f(2+s,s) - |s|^{p-1}  \frac{4}{p-1} f(2+s,s)\ ds $$
which by the product rule is equal to
\begin{equation}\label{vatr}
 V(t,r) = \frac{1+r}{2} f(t,r) + \frac{1}{2|r|^{p-1}} \int_r^0 |s|^{p-1} \left[ \frac{p-5}{p-1} + \frac{(p-1)(1+s)}{s} \right] f(2+s,s)\ ds.
\end{equation}
Note that $f(2+s,s)$ is only non-zero when $s = -1+O(\delta)$, in which case it is of size $O( A \delta^{-\frac{4}{p-1}} )$ thanks to \eqref{eb}, \eqref{ec}.  This gives the upper bound in \eqref{dipp}.  Now we turn to the lower bound.  First suppose that $-(1-\frac{\delta}{2})t \leq r$, then $f$ is non-negative in all of its appearances in \eqref{vatr}.  As we are in the supercritical case $p>5$, the factor $\frac{p-5}{p-1} + \frac{(p-1)(1+s)}{s}$ is positive (indeed it is $\gtrsim 1$) for $\delta$ small enough, and the claim follows in this case.

It remains to consider the case when $-t \leq r \leq -(1-\frac{\delta}{2}) t$.  In this case we can use the lower bound
$$ f(t,r) \geq - (t+r)^{-\frac{4}{\delta}}$$
and conclude that the term $\frac{1+r}{2} f(t,r)$ is at least $-O( \delta^{\frac{p-5}{p-1}} )$.  A similar argument shows that the contribution to \eqref{vatr} coming from those $s$ with $-(2+s) \leq s \leq -(1-\frac{\delta}{2}) (2+s)$ is at least $-O( \delta^{\frac{p-5}{p-1}} )$.  On the other hand, from \eqref{big}, the contribution of those $s$ with $s >-(1-\frac{\delta}{2}) (2+s)$ is $\gtrsim (A-O(1)) \delta^{\frac{p-5}{p-1}}$.  As $A$ is assumed to be large, the claim follows.
\end{proof}

On the support of $V$ in $\Gamma_1^l$, we see from \eqref{ed}, \eqref{etr-def} that
$$ 
(\partial_t - \partial_r) e_+ = O( A t^{-\frac{p+3}{p-1}} \delta^{-\frac{4}{p-1}} )
$$
and hence by \eqref{storp} and Proposition \ref{vpos}
\begin{equation}\label{vanity}
(\partial_t + \partial_r) V = O( A t^{-\frac{p+3}{p-1}} \delta^{-\frac{4}{p-1}} ).
\end{equation}

\subsection{Selection of $V$ in the right half of the cone}  

Once $V$ has been constructed in the left half $\Gamma^l_1$ of the light cone, we extend it to the right half $\Gamma^r_1 := \{ (t,r) \in \Gamma_1: r \geq 0 \}$ by even extension, thus
$$ V(t,r) := V(t,-r)$$
for all $(t,r) \in \Gamma^r_1$.  Since $V$ vanished for $-(1-\delta) t \leq r \leq 0$, we see that $V$ is smooth on all of $\Gamma_1$, and vanishing in the interior cone $\{ (t,r) \in \Gamma_1: |r| \leq (1-\delta) t \}$.  It also obeys the non-negativity property \eqref{rpos}.  From reflecting \eqref{vanity} and Proposition \ref{vpos} we have the bounds
\begin{equation}\label{vpos-2}
V = O( A t^{-\frac{4}{p-1}} \delta^{\frac{p-5}{p-1}} )
\end{equation}
and
\begin{equation}\label{vanity-2}
(\partial_t - \partial_r) V = O( A t^{-\frac{p+3}{p-1}} \delta^{-\frac{4}{p-1}} )
\end{equation}
when $(1-\delta) t \leq r \leq t$.

\subsection{Selection of $e_+$ in the right half of the cone}

Thus far, $V$ has been defined on all of $\Gamma_1$, and $e_+$ defined on $\Gamma_1^l$.  We now extend $e_+$ to $\Gamma_1^r$ by solving \eqref{storp}, or more precisely by setting
\begin{equation}\label{votr}
 e_+(t,r) := e_+(t+r,0) + \int_0^r ((\partial_t + \partial_r) V)(t+r-s,s) + \frac{p-1}{s} V(t+r-s, s)\ ds
\end{equation}
for $0 < r \leq t$;
note that the integral is well-defined since $V$ vanishes near the time axis.  One easily checks that $e_+(t,r) = (t+r)^{-\frac{4}{p-1}}$ for $0 \leq r \leq (1-\delta) t$, and so $e_+$ extends smoothly to all of $\Gamma_1$ and is equal to $(t+r)^{-\frac{4}{p-1}}$ in the interior cone $\{ (t,r) \in \Gamma_1: |r| \leq (1-\delta) t \}$.  It is also clear from construction that $e_+$ is homogeneous of order $-\frac{4}{p-1}$.  From the fundamental theorem of calculus we see that $e_+$ and $V$ obey \eqref{storp} on $\Gamma_1^r$, and hence on all of $\Gamma_1$.  From \eqref{vpos-2}, \eqref{vanity-2} we see that the integrand is of size $O( A t^{-\frac{p+3}{p-1}} \delta^{-\frac{4}{p-1}} )$ when $r = (1-O(\delta)) t$, and vanishes otherwise, which leads (for $\delta$ small enough) to the crude upper and lower bounds
\begin{equation}\label{cruel}
t^{-\frac{p+3}{p-1}} \lesssim e_+(t,r) \lesssim t^{-\frac{p+3}{p-1}}
\end{equation}
throughout $\Gamma_1^r$.

\subsection{Selection of $e_-$ and $\tilde E_{tr}$}

We reflect the function $e_+$ around the time axis to create a new function $e_-\colon \Gamma_1 \to \R$:
$$ e_-(t,r) := e_+(t,r).$$
Like $e_+$, the function $e_-$ is smooth and homogeneous of order $-\frac{4}{p-1}$.  It equals $(t-r)^{-\frac{4}{p-1}}$ in the interior cone $\{ (t,r) \in \Gamma_1: |r| \leq (1-\delta) t \}$.  On $\Gamma_1^l$ it obeys the crude upper and lower bounds
\begin{equation}\label{cruel-refl}
t^{-\frac{p+3}{p-1}} \lesssim e_-(t,r) \lesssim t^{-\frac{p+3}{p-1}}
\end{equation}
and in the region $(1-\delta) t \leq r \leq t$ we have the bounds
\begin{equation}\label{jam}
(\delta t)^{-\frac{4}{p-1}} \lesssim e_-(t,r) \lesssim A (\delta t)^{-\frac{4}{p-1}} 
\end{equation}
thanks to \eqref{ec}.

Recall from \eqref{null-energy} that the field $e_+$ is intended to ultimately be of the form $\frac{1}{2} \tilde E_{tt} + \frac{1}{2} \tilde E_{rr} + \tilde E_{tr}$.  Similarly, $e_-$ is intended to be of the form 
\begin{equation}\label{intend}
e_- = \frac{1}{2} \tilde E_{tt} + \frac{1}{2} \tilde E_{rr} - \tilde E_{tr}.
\end{equation}
Accordingly, we may now define $\tilde E_{tr}$ as
\begin{equation}\label{rrtt}
 \tilde E_{tr} := \frac{e_+-e_-}{2}.
\end{equation}
This is clearly smooth, odd, and homogeneous of order $-\frac{4}{p-1}$.
We also see that the quantity $\tilde E_{tt} + \tilde E_{rr}$ is now specified:
\begin{equation}\label{ttrr}
\tilde E_{tt} + \tilde E_{rr} = e_+ + e_-.
\end{equation}

We are left with two remaining unknown scalar fields to specify: the mass density $\tilde M$ and the energy equipartition $-\tilde E_{tt} + \tilde E_{rr}$, which determines the fields $\tilde E_{tt}$ and $\tilde E_{rr}$ by \eqref{ttrr}.  The requirements needed for Theorem \ref{main-6} that have not already been verified are as follows:

\begin{itemize}
\item $\tilde M$ is smooth, even, and homogeneous of order $\frac{2p-6}{p-1}$; $-\tilde E_{tt}+\tilde E_{rr}$ is smooth, even, and homogeneous of order $-\frac{4}{p-1}$.
\item $\tilde M, \tilde E_{tt}$ vanishes to second order at $r=0$, $\tilde E_{tr} - \frac{1}{2r} \partial_t \tilde M$ vanishes to third order, and $\tilde E_{rr} - \frac{1}{r} \partial_r \tilde M + \frac{1}{r^2} \tilde M$ to fourth order.  
\item One has the equations \eqref{vdef} and \eqref{null-energy} (and hence also \eqref{intend}).
\item The matrix \eqref{la-4} is strictly positive definite for $r \neq 0$, and for $r=0$ one has the condition \eqref{mur}.
\end{itemize}

As there is only one equation (beyond homogeneity and reflection symmetry) constraining $\tilde M$ and $-\tilde E_{tt}+\tilde E_{rr}$ - namely, \eqref{vdef} - the problem of selecting these two fields is underdetermined, and thus subject to a certain amount of arbitrary choices.  We will select these fields first in the exterior region $\{ (t,r) \in \Gamma: |r| \geq t/2 \}$, and then fill in the interior using a different method.

\subsection{Selection of $M, -\tilde E_{tt}+ \tilde E_{rr}$ away from the time axis}

In the exterior region $\{ (t,r) \in \Gamma: |r| \geq t/2 \}$, we shall simply select the field $\tilde M$ to be a small but otherwise rather arbitrary field, and then use \eqref{vdef} to determine $-\tilde E_{tt}+ \tilde E_{rr}$.

More precisely, let $\tilde M(1,r)$ be a smooth even function on the region $\{ r: 1/2 \leq |r| \leq 1 \}$ obeying the following properties:
\begin{itemize}
\item For $1/2 \leq |r| \leq 3/4$, one has
\begin{equation}\label{mir}
 \tilde M(1,r) = \delta ( (1+r)^{\frac{2p-6}{p-1}} + (1-r)^{\frac{2p-6}{p-1}} ).
\end{equation}
(This condition will not be used directly in this part of the construction, but is needed for compatibility with the next part.)
\item For $1/2 \leq |r| \leq 1$, one has the bounds
\begin{equation}\label{mir-1}
\delta \lesssim \tilde M(1,r) \lesssim \delta
\end{equation}
and
\begin{equation}\label{mir-2}
\frac{d}{dr} \tilde M(1,r), \frac{d^2}{dr^2} \tilde M(1,r) = O(\delta).
\end{equation}
\end{itemize}
It is clear that one can select such a function.  We then extend $\tilde M$ to $\{ (t,r) \in \Gamma_1: |r| \geq t/2 \}$ by requiring that $\tilde M$ be homogeneous of order $\frac{2p-6}{p-1}$.  Then $\tilde M$ is smooth and even, and one has the bounds
\begin{align}
\delta t^{\frac{2p-6}{p-1}} \lesssim \tilde M(t,r) &\lesssim \delta t^{\frac{2p-6}{p-1}} \label{mrt-1} \\
\frac{d}{dr} \tilde M(t,r), \frac{d}{dt} \tilde M(t,r) &= O( \delta t^{\frac{p-5}{p-1}} ) \label{mrt-2} \\
\frac{d^2}{dr^2} \tilde M(t,r), \frac{d^2}{dt^2} \tilde M(t,r) &= O( \delta t^{-\frac{4}{p-1}} ) \label{mrt-3}
\end{align}
in the region $\{ (t,r) \in \Gamma_1: |r| \geq t/2 \}$.  

We then define $-\tilde E_{tt} + \tilde E_{rr}$ on this region by enforcing \eqref{vdef}, thus
\begin{equation}\label{etra}
-\tilde E_{tt} + \tilde E_{rr} := \frac{1}{2} ( - \partial_{tt} \tilde M + \partial_{rr} \tilde M ) - (p+1) V.
\end{equation}
Combining this with \eqref{ttrr}, this defines $\tilde E_{tt}$ and $\tilde E_{rr}$.  It is easy to see that these fields are smooth, even and homogeneous of order $-\frac{4}{p-1}$ on $\{ (t,r) \in \Gamma_1: |r| \geq t/2 \}$.

We now claim that the matrix \eqref{la-4} is strictly positive definite in the region $\{ (t,r) \in \Gamma_1: |r| \geq t/2 \}$.  By homogeneity and reflection symmetry, it suffices to verify this when $t=1$ and $1/2 \leq r \leq 1$.  Using the identity
$$
\begin{pmatrix}
\tilde M & \frac{1}{2} (\partial_t+\partial_r) \tilde M & \frac{1}{2} (\partial_t-\partial_r) \tilde M \\
\frac{1}{2} (\partial_t+\partial_r) \tilde M & 2e_+ & -\tilde E_{tt} + \tilde E_{rr} \\
\frac{1}{2} (\partial_t-\partial_r) \tilde M & -\tilde E_{tt} + \tilde E_{rr} & 2e_-
\end{pmatrix}
= 
\begin{pmatrix}
1 & 0 & 0 \\
0 & 1 & 1 \\
0 & -1 & 1
\end{pmatrix}
\begin{pmatrix}
\tilde M & \frac{1}{2}  \partial_t \tilde M & \frac{1}{2} \partial_r \tilde M \\
\frac{1}{2} \partial_t \tilde M & \tilde E_{tt} & \tilde E_{tr} \\
\frac{1}{2} \partial_r M & \tilde E_{tr} & \tilde E_{rr}
\end{pmatrix}
\begin{pmatrix}
1 & 0 & 0 \\
0 & 1 & -1 \\
0 & 1 & 1
\end{pmatrix},
$$
it suffices to show that the matrix
$$
\begin{pmatrix}
\tilde M & \frac{1}{2} (\partial_t+\partial_r) \tilde M & \frac{1}{2} (\partial_t-\partial_r) \tilde M \\
\frac{1}{2} (\partial_t+\partial_r) \tilde M & 2e_+ & -\tilde E_{tt} + \tilde E_{rr} \\
\frac{1}{2} (\partial_t-\partial_r) \tilde M & -\tilde E_{tt} + \tilde E_{rr} & 2e_-
\end{pmatrix}
$$
is strictly positive definite.

If $r \leq 1-\delta$, then all off-diagonal terms are $O(\delta)$ thanks to \eqref{mir-2}, \eqref{etra}, while the diagonal terms are $\gtrsim \delta$, $\gtrsim 1$, and $\gtrsim 1$ respectively, and the positive definiteness is easily verified, since the associated quadratic form is at least
$$ \gtrsim \delta x_1^2 + x_2^2 + x_3^2 - O( \delta |x_1| |x_2| ) - O( \delta |x_1| |x_3| ) - O( \delta |x_2| |x_3| )$$
which is easily seen to be positive for $\delta$ small enough.  If $r < 1-\delta$, then the off-diagonal terms are $O(\delta)$ in the top row and left column, and $O( A \delta^{\frac{p-5}{p-1}} )$ in the bottom right minor by \eqref{vpos-2}, while the diagonal terms are $\gtrsim \delta$, $\gtrsim 1$, and $\gtrsim \delta^{-\frac{4}{p-1}}$ by \eqref{mir-1}, \eqref{cruel}, \eqref{jam}, so the associated quadratic form is
$$ 
\gtrsim \delta x_1^2 + x_2^2 + \delta^{-\frac{4}{p-1}} x_3^2 - O( \delta |x_1| |x_2| ) - O( \delta |x_1| |x_3| ) - O( A \delta^{\frac{p-5}{p-1}} |x_2| |x_3| )$$
which is again positive definite (note that $A\delta^{\frac{p-5}{p-1}}$ can be chosen to be much smaller than the geometric mean of $\delta$ and $\delta^{-\frac{4}{p-1}}$).

\subsection{Selection of $M, \tilde E_{tt}, \tilde E_{rr}$ near the time axis}

Now we restrict attention to the interior region $\Gamma_1^i := \{ (t,r) \in \Gamma_1: |r| \leq t/2 \}$; all identities and estimates here are understood to be on this region unless otherwise specified.

We will now reverse the Gram matrix reduction from previous sections, and construct $\tilde M, \tilde E_{tt}, \tilde E_{rr}$ in $\Gamma_1^i$ from an (infinite-dimensional) vector-valued solution to the (free, $1+1$-dimensional) wave equation.  Let $H$ be a Hilbert space and let $t \mapsto f(t)$ be a family of vectors $f(t)$ in $H$ smoothly parmeterised by a parameter $t \in (0,+\infty)$ (so that all derivatives in $t$ exist in the strong sense and are continuous); we will select this family more precisely later.  We introduce the smooth vector-valued field $\phi\colon \Gamma_1^i \to H$ by the formula
$$ \phi(t,r) := f(t+r) - f(t-r) $$
and we will define $\tilde M, \tilde E_{tt}, \tilde E_{rr}\colon \Gamma_1^i \to \R$ by the formulae
\begin{align*}
\tilde M(t,r) &:= \langle \phi(t,r), \phi(t,r) \rangle_H \\
\tilde E_{tt}(t,r) &:= \langle \partial_t \phi(t,r), \partial_t \phi(t,r) \rangle_H \\
\tilde E_{rr}(t,r) &:= \langle \partial_r \phi(t,r), \partial_r \phi(t,r) \rangle_H.
\end{align*}
Since $\phi$ is smooth and odd in $r$, these functions are smooth and even in $r$.  If we impose the additional hypothesis that the Gram matrix $\langle f(s), f(t) \rangle_H$ has the scaling symmetry
\begin{equation}\label{fast}
 \langle f(\lambda s), f(\lambda t) \rangle_H = \lambda^{\frac{2p-6}{p-1}} \langle f(s), f(t) \rangle_H
\end{equation}
for $s,t,\lambda > 0$, then $M$ will be homogeneous of order $\frac{2p-6}{p-1}$; furthermore, by differentiating \eqref{fast} with respect to both $s$ and $t$ we see that
\begin{equation}\label{fast2}
 \langle f'(\lambda s), f'(\lambda t) \rangle_H = \lambda^{-\frac{4}{p-1}} \langle f'(s), f'(t) \rangle_H
\end{equation}
(where $f'$ denotes the derivative of $f$) and so $\tilde E_{tt}, \tilde E_{rr}$ will be homogeneous of order $-\frac{4}{p-1}$.

Observe that
\begin{align*}
\frac{1}{2} \tilde E_{tt} + \frac{1}{2} \tilde E_{rr} + \langle \partial_t \phi, \partial_r \phi \rangle_H
&= \frac{1}{2} \| (\partial_t + \partial_r) \phi \|_H^2 \\
&= 2 \| f'( t+r ) \|_H^2
\end{align*}
and similarly
$$ \frac{1}{2} \tilde E_{tt} + \frac{1}{2} \tilde E_{rr} - \langle \partial_t \phi, \partial_r \phi \rangle_H = 2 \|f'(t-r)\|_H^2.$$
Thus, if we impose the additional normalisation
\begin{equation}\label{hypo}
\| f'(1) \|_H = \frac{1}{\sqrt{2}} 
\end{equation}
and hence by \eqref{fast2}
\begin{equation}\label{stew}
\| f'(t) \|_H = \frac{1}{\sqrt{2}} t^{-\frac{2}{p-1}} 
\end{equation}
we see from the identities $e_\pm(t,r) = (t \pm r)^{-\frac{4}{p-1}}$ in $\Gamma_1^i$ that
$$ \frac{1}{2} \tilde E_{tt} + \frac{1}{2} \tilde E_{rr} \pm \langle \partial_t \phi, \partial_r \phi \rangle_H = e_\pm.$$
In particular, \eqref{ttrr} holds, and from \eqref{rrtt} one has
$$ \tilde E_{tr} = \langle \partial_t \phi, \partial_r \phi \rangle_H.$$
We also obtain the equations \eqref{null-energy} and \eqref{intend}.

Next, it is clear that $\phi$ solves the wave equation
$$ - \partial_{tt} \phi + \partial_{rr} \phi = 0$$
so in particular
$$ \langle \phi, - \partial_{tt} \phi + \partial_{rr} \phi  \rangle_H = 0 $$
which implies in particular (cf. \eqref{pr}) that
$$ \frac{1}{2} ( - \partial_{tt} \tilde M + \partial_{rr} \tilde M ) + \tilde E_{tt} - \tilde E_{rr} = 0.$$
Since $V$ vanishes on $\Gamma_1$, we conclude that \eqref{vdef} holds.

Next, from differentiating the formula for $\tilde M$, one has
$$ \frac{1}{2} \partial_t \tilde M = \langle \phi, \partial_t \phi \rangle_H$$
and
$$ \frac{1}{2} \partial_r \tilde M = \langle \phi, \partial_t \phi \rangle_H$$
and so the quadratic form associated with \eqref{la-4} factorises as
$$ \| x_1 \phi + x_2 \partial_t \phi + x_3 \partial_r \phi \|_H^2.$$
This is clearly positive semi-definite at least; to make it positive definite for $r \neq 0$, it will suffice to enforce the condition
\begin{equation}\label{linind}
f(s), f(t), f'(s), f'(t) \hbox{ linearly independent}
\end{equation}
for all distinct $s,t > 0$.

Suppose we assume the long-range orthogonality condition
\begin{equation}\label{long}
\langle f(s), f(t) \rangle_H = 0
\end{equation}
whenever $t/s > 1.1$ or $s/t > 1.1$.   Then in the region $\{ (t,r) \in \Gamma_1^i: |r| \geq t/4 \}$ away from the time axis, we have from Pythagoras' theorem that
$$ 
\tilde M(t,r) = \| f(t+r) \|_H^2 + \|f(t-r) \|_H^2.$$
In particular, if we also impose the normalisation
\begin{equation}\label{fail}
\|f(1)\|_H = \sqrt{\delta}
\end{equation}
then (from \eqref{fast}) we have
$$ \tilde M(t,r) = \delta ( (t+r)^{\frac{2p-6}{p-1}} + (t-r)^{\frac{2p-6}{p-1}} ) $$
in the region $\{ (t,r) \in \Gamma_1^i: |r| \geq t/4 \}$.  In particular from \eqref{mir} and homogeneity we see that $\tilde M$ on $\Gamma_1^i$ joins up smoothly with its counterpart in the exterior region $\{ (t,r) \in \Gamma_1: |r| \geq t/2 \}$; by \eqref{vdef} we see that $-\tilde E_{tt}+\tilde E_{rr}$ does too. By \eqref{ttrr}, \eqref{rrtt} we now see that all of the fields $M, \tilde E_{tt}, \tilde E_{rr}, \tilde E_{tr}$ are smooth on all of $\Gamma_1$.

Now we study the vanishing properties of the various fields constructed at $r=0$, for a fixed value of $t$.  From Taylor expansion we have
$$ \phi(t,r) = 2r f'(t) + \frac{1}{3} r^3 f'''(t) + O( |r|^5 ) $$
as $r \to 0$ (where the error term denotes a quantity in $H$ of norm $O(|r|^5)$, and the implied constant can depend on $t$ and $\phi$.  Furthermore, these asymptotics behave in the expected fashion with respect to differentiation in time or space, thus for instane
$$ \partial_r \phi(t,r) = 2 f'(t) + r^2 f'''(t) + O( |r|^4 ) $$
and
$$ \partial_t \phi(t,r) = 2r f''(t) + \frac{1}{3} r^3 f^{(4)}(t) + O( |r|^5 ).$$
Taking inner products, we conclude the asymptotics
\begin{align*}
 \tilde M(t,r) &= 4r^2 \|f'(t)\|_H^2 + \frac{4}{3} r^4 \langle f'(t), f'''(t) \rangle_H + O( |r|^6 ) \\
 \tilde E_{tt}(t,r) &= 4r^2 \|f''(t)\|_H^2 + O( |r|^4 ) \\
 \tilde E_{tr}(t,r) &= 4r \langle f'(t), f''(t) \rangle_H + O( |r|^3 )\\
 \tilde E_{rr} &= 4 \|f'(t)\|_H^2 + 4 r^2 \langle f'(t), f'''(t) \rangle_H + O( |r|^4 ).
\end{align*}
The asymptotic for $\tilde M$ behaves well with respect to derivaties, thus for instance
$$
\partial_t \tilde M(t,r) = 8 r^2 \langle f'(t), f''(t) \rangle_H + O(|r|^4)
$$
and
$$ 
\partial_r \tilde M(t,r) = 8 r \|f'(t)\|_H^2 + \frac{16}{3} r^3 \langle f'(t), f'''(t) \rangle_H + O(|r|^5).$$
Among other things, this shows (using \eqref{hypo}) that the condition \eqref{mur} reduces to
\begin{equation}\label{mur-2}
 \|f''(1)\|_H > \frac{2}{p-1} \frac{1}{\sqrt{2}}.
\end{equation}
It is also clear from these asymptotics that $\tilde E$ and $\tilde E_{tt}$ vanish to second order, and $\tilde E_{tr} - \frac{1}{2r} \partial_t \tilde M$ vanishes to third order; a brief calculation also shows that $\tilde E_{rr} - \frac{1}{r} \partial_r \tilde M + \frac{1}{r^2} \tilde M$ vanishes to fourth order.

To summarise: in order to conclude all the required properties for Theorem \ref{main-6}, it suffices to locate a smooth curve $t \mapsto f(t)$ in a Hilbert space $H$ which obeys the hypotheses \eqref{fast}, \eqref{hypo}, \eqref{linind}, \eqref{long}, \eqref{fail}, \eqref{mur-2}.

We take the Hilbert space $H$ to be the space $L^2(\R)$ of square-integrable real-valued functions on $\R$ with Lebesgue measure.  The functions $f(t) \in H$ will take the form
$$ f(t)(x) := t^{\frac{p-3}{p-1}} \psi( x - \log t )$$
where $\psi\colon \R \to \R$ is a bump function whose (closed) support is precisely $[0,0.01]$ (that is to say, the set $\{ \psi \neq 0\}$ is a dense subset of $[0,0.01]$) depending on $\delta$ and $p$ to be chosen shortly.  It is clear from construction that \eqref{fast} and \eqref{long} hold.  The condition \eqref{fail} becomes
$$
\int_\R \psi(x)^2\ dx = \delta $$
while the condition \eqref{hypo} becomes
$$
\int_\R \psi'(x)^2\ dx = \frac{1}{2}.
$$
It is easy to see that we can select $\psi$ with closed support precisely $[0,0.01]$ with both of these normalisations, basically because the Dirichlet form $\langle \phi', \psi' \rangle$ is unbounded on $L^2([0,0.01])$.

Now we verify the linear independence claim \eqref{linind}.  We may assume without loss of generality that $s=1$ and $t>1$.  Then we have a linear dependence between $\psi$ and $\psi'$ in a neighbourhood of $0$; since $\psi, \psi'$ vanish to the left of $0$, the Picard uniqueness theorem for ODEs then implies that $\psi$ vanishes a little to the right of $0$ also, contradicting the hypothesis that $\psi$ has closed support containing $0$.  This gives \eqref{linind}.

A similar argument shows that $f'(1)$ amd $f''(1)$ are linearly independent.  Squaring and differentiating \eqref{stew} at $t=1$ gives
$$ \langle f'(1), f''(1) \rangle_H = - \frac{2}{p-1} \frac{1}{2} $$
and \eqref{mur-2} then follows from \eqref{hypo} and the Cauchy-Schwarz inequality, using the linear independence to get the strict inequality.  This (finally) completes the proof of Theorem \ref{main-6} and hence Theorem \ref{main}.

\end{document}